\documentclass[12pt]{amsart}

\usepackage{amsmath,amssymb,epic,lscape}

\newtheorem{theorem}{Theorem}[section]

\theoremstyle{definition}
\newtheorem{definition}[theorem]{Definition}
\newtheorem{example}[theorem]{Example}

\theoremstyle{remark}
\newtheorem{remark}[theorem]{Remark}

\numberwithin{equation}{section}
\providecommand{\bysame}{\leavevmode\hbox to3em{\hrulefill}\thinspace}

\begin{document}

\def\DJ{{\hbox{D\kern-.8em\raise.15ex\hbox{--}\kern.35em}}}
\def\DJo{$\;$\kern-.4em
    \hbox{D\kern-.8em\raise.15ex\hbox{--}\kern.35em okovi\'c}}

\def\NSERC{Supported in part by the NSERC Grant A-5285.}

\font\germ=eufm10
\def\sL{{\mbox{\germ sl}}}

\def\al{{\alpha}}
\def\be{{\beta}}
\def\vf{{\varphi}}
\def\la{{\lambda}}
\def\bR{{\mbox{\bf R}}}
\def\bZ{{\mbox{\bf Z}}}
\def\bC{{\mbox{\bf C}}}
\def\bN{{\mbox{\bf N}}}
\def\bH{{\mbox{\bf H}}}
\def\pA{{\mathcal A}}
\def\pF{{\mathcal F}}
\def\pE{{\mathcal E}}
\def\pC{{\mathcal C}}
\def\pX{{\mathcal X}}
\def\pY{{\mathcal Y}}
\def\pZ{{\mathcal Z}}
\def\Cb{{\bar{C}}}
\def\tr{{\rm tr\;}}
\def\Tr{{\rm Tr\;}}
\def\Aut{{\mbox{\rm Aut}}}
\def\Sk{{\mbox{\rm Skew}}}
\def\GL{{\mbox{\rm GL}}}
\def\SL{{\mbox{\rm SL}}}
\def\SO{{\mbox{\rm SO}}}
\def\Sp{{\mbox{\rm Sp}}}
\def\Un{{\mbox{\rm U}}}
\def\Ort{{\mbox{\rm O}}}

\renewcommand{\subjclassname}{\textup{2000} Mathematics Subject
Classification }

\title[Cyclic difference families with two base blocks]
{ Cyclic ${ \mathbf {(v;r,s;\la)} }$ difference families with
two base blocks and ${ \mathbf {v\le50} }$ }

\author[D.\v{Z}. \DJ okovi\'{c}]
{Dragomir \v{Z}. \DJ okovi\'{c}}

\address{Department of Pure Mathematics, University of Waterloo,
Waterloo, Ontario, N2L 3G1, Canada}

\email{djokovic@uwaterloo.ca}

\thanks{
The author was supported by an NSERC Discovery Grant.}

\keywords{Difference family, supplementary difference sets,
balanced incomplete block designs, genetic algorithm}

\date{}

\begin{abstract}
We construct many new cyclic $(v;r,s;\la)$ difference families
with $v\ge2r\ge2s\ge4$ and $v\le50$. In particular we
construct the difference families with parameters
\begin{center}
\begin{tabular}{lll}
(45;18,10;9),& (45;22,22;21),& (47;21,12;12), \\
(47;19,15;12),& (47;22,14;14),& (48;20,10;10), \\
(48;24,4;12), & (50;25,20;20) &
\end{tabular}
\end{center}
for which the existence question was an open problem.

We point out that the $(45;22,22;21)$ difference family
gives a balanced incomplete block design (BIBD)
with parameters $v=45$, $b=90$, $r=44$, $k=22$ and $\la=21$,
and that the one with parameters (50;25,20;20) gives a pair
of binary sequences of length 50 with zero periodic
autocorrelation function (the periodic analog of a Golay pair).
The new SDSs include nine new D-optimal designs. 

A normal form for cyclic difference families (with base blocks of
arbitrary sizes) is proposed and used effectively in compiling
our selective listings in Tables 3-6 of known and new
difference families in the above range.
\end{abstract}

\maketitle
\subjclassname{ 05B20, 05B30 }
\vskip5mm

\section{Introduction}

We consider difference families in finite abelian groups
whose base blocks may be of different sizes,
also known as supplementary difference sets (SDS).
A few infinite families have been constructed,
most of them by using cyclotomy.
They have been studied for long time and have been used
to construct balanced incomplete block designs (BIBD),
Hadamard matrices, skew-Hadamard matrices, and other designs.
We recommend \cite{HCD,GS1,GS2,SY,RW} for an overview of this topic.

We shall restrict our scope here to the case of SDSs with
exactly two base blocks. Their sizes will be denoted by $r$ and $s$.
The underlying group will be cyclic of order $v$, identified
with $\bZ_v$. As usual, we attach to such SDS its
parameters $(v;r,s;\la)$ and the order $n$,
where $\la$ is the index of the family and $n=r+s-\la$.
For precise definition see the next section.
The systematic search for these SDSs in the range
$v\le50$, but restricted to odd values of $v$, was initiated
in \cite{CCS}. In our paper \cite{DZ1} we have continued this
work and extended the search to include the even values of $v$.

There is a simple necessary condition on $(v;r,s;\la)$ for the
existence of an SDS with this set of parameters (see (\ref{jedn})
below). We refer to the parameters satisying that condition
as feasible parameters. In a more recent joint paper \cite{MDV}
the question of existence of SDSs with feasible parameter sets
has been settled for all $v\le40$. In the range
$40<v\le50$ there remained 23 undecided cases.

Our objective in this paper is twofold. First we feel that there
is a need to collect the known results in ready for use
tabular form. For that purpose we introduce a normal form 
for cyclic difference families (with arbitrary number of base
blocks of various sizes) in order to be able to compare the
results from different sources and avoid the duplication.
While doing this we also constructed many new SDSs, not
equivalent to the known ones. For instance, 
each of the Tables 5 and 6 contains 59 representatives 
of new equivalence classes of SDSs. These two tables contain
new D-optimal designs with parameters

\begin{center}
\begin{tabular}{lll}
(31;15,10;10),& (37;16,13;11),& (41;16,16;12), \\ 
(43;18,16;13),& (43;21,15;15), & (49;22,18;16)
\end{tabular}
\end{center}
(nine new designs in total).

Our second objective is to make further contribution to the question
of existence of SDSs with the given feasible set of parameters.
We have constructed SDSs proving the existence in the following
8 of the 23 previously undecided cases with $v\le50$:

\begin{tabular}{llll}
&&& \\
(45;18,10;9), & (45;22,22;21), & (47;21,12;12), & (47;19,15;12), \\
(47;22,14;14), & (48;20,10;10), & (48;24,4;12), & (50;25,20;20). \\
&&& \\
\end{tabular}

We pay special attention to three types of these difference
families that were investigated separately in the past.
The first type are the SDSs with $r=s$, i.e., the difference families
with two base blocks of the same size. They are important because
they can be used to construct BIBD's. For instance, the
(45;22,22;21) difference family gives rise to a BIBD with
parameters mentioned in the abstract. The second type are the SDSs
with $v=2n+1$. They are used to construct D-optimal designs
of order $2v$ of circulant i.e., Ehlich type. The third type
are the SDSs with $v=2n$. They are equivalent to pairs of binary
sequences of length $v$ having zero periodic autocorrelation
function (see \cite{BA} for a more general fact).
The latter are important in various engineering
applications. Apart from the lengths $v$ of known Golay pairs,
there are only two known even integers $v$ for which such
binary sequences are known: 34 and 50. The ones with $v=50$
are constructed in this paper.

Our results are presented in Tables 3-6. The most
interesting ones are in the last table.

All new SDSs recorded in Tables 3-6 have been constructed by
running our genetic type algorithm, sometimes for several days.
For a short description of the algorithm we refer the reader
to \cite{MDV}, but see also \cite{DA}.

\section{Equivalence of difference families}

As we are interested in recording known cyclic difference families
having exactly two base blocks (not necessarily distinct or of the
same size), we are faced
with the problem of testing different families for equivalence.
Of course, first we have to define the notion of equivalence
for difference families.
Then, we have to construct a suitable normal form.
We shall be more general and allow any number of base blocks.

Let $v$ be a positive integer and denote by $\bZ_v=\{0,1,\ldots,v-1\}$
the ring (and also the additive group) of integers modulo $v$.
By $\bZ_v^*$ we denote the (multiplicative)
group of invertible elements of $\bZ_v$.
To any subset $X$ of $\bZ_v$ we associate a function
$N_X:\bZ_v\to\bZ$ whose value at a point $a\in\bZ_v$ is equal to
the cardinal of the set $\{(x,y)\in X\times X: y-x=a \}$.
Note that $N_X(0)=|X|$, the size of $X$.

\begin{definition} \label{df}
We say that a sequence $\pX=(X_1,\ldots,X_m)$, with $X_i\subseteq\bZ_v$
for all $i$, is a \emph{difference family} or that $X_1,\ldots,X_m$
are \emph{supplementary difference sets} (SDS) if the function
\[ N_\pX=N_{X_1}+\cdots+N_{X_m} \]
takes the same value, say $\la$, on all nonzero elements $a\in\bZ_v$.
In that case, if $k_i=|X_i|$, we say that
\begin{equation} \label{param}
(v;k_1,\ldots,k_m;\la)
\end{equation}
is the \emph{set of parameters} of $\pX$, the $X_i$'s are its
\emph{base blocks}, $\la$ is its \emph{index} and
\[ n=k_1+\cdots+k_m-\la \]
is its \emph{order}.
If also $m=1$, we say that $X_1$ is a \emph{difference set}.
We denote by $\pF_m$ the set of all difference families
$(X_1,\ldots,X_m)$ in $\bZ_v$.
\end{definition}

Let us introduce the elementary transformations 
$\pX\to\pY$ on the set $\pF_m$. There are five types
of such transformations.
For $X\subseteq\bZ_v$, we denote by $X^\dag$ the symmetric
difference of $X$ and the set of odd integers in $\bZ_v$.

Let $\pX=(X_1,\ldots,X_m)\in\pF_m$ and let $\pY=(Y_1,\ldots,Y_m)$
be obtained from $\pX$ by one of the following
\emph{elementary transformations}: 
\begin{itemize}
\item[(i)] Replace an $X_i$ with $X_i+t$ for some $t\in\bZ_v$.
\item[(ii)] For some $a\in\bZ_v^*$, replace each $X_i$ with $a\cdot X_i$.
\item[(iii)] Replace an $X_i$ with $(-1)\cdot X_i$.
\item[(iv)] Replace an $X_i$ with its complement in $\bZ_v$.
\item[(v)] If $v$ is even and $mv=4n$, replace each $X_i$ with $X_i^\dag$.
\end{itemize}
It is easy to verify that $\pY\in\pF_m$. If needed, consult
e.g. \cite{BA} for more details about the case (v).

\begin{definition} \label{equiv}
We say that $\pX,\pY\in\pF_m$ are \emph{equivalent} if there exists
a finite sequence of elementary transformations which sends
$\pX$ to $\pY$. This is an equivalence relation on $\pF_m$.
\end{definition}

Note that equivalent SDSs $\pX$ and $\pY$ may have different
sets of parameters. For instance, if $\pX$ has parameters
(\ref{param}) and $\pY$ is obtained from $\pX$ by the elementary
transformation (iv), then the set of parameters of $\pY$ is
\[ (v;k_1,\ldots,k_{i-1},v-k_i,k_{i+1},\ldots,k_m;\la+v-2k_i). \]

In the case when $m=2$ and $v$ is odd, this definition of equivalence
is different from the one adopted in \cite{KK1},
which additionally allows the $X_i$'s to be permuted.

\section{Normal form for difference families}

Let $X$ be a $k$-subset of $\bZ_v$ and let $d$ be any nonnegative
integer. As we have identified
$\bZ_v$ with a subset of $\bZ$, it makes sense to raise each
$x\in X$ to power $d$ in the ring $\bZ$ and add all these
powers in $\bZ$. By convetion, $x^0=1$ for all $x\in\bZ_v$.
We denote this sum by $\sigma_d(X)$. Thus
\[ \sigma_d(X)=\sum_{x\in X} x^d \in\bZ. \]

Now let $X$ and $Y$ be any two subsets of $\bZ_v$. If
$\sigma_d(X)=\sigma_d(Y)$ for all $0\le d<v$, then
$X=Y$. We shall write $X<Y$ if there exists a nonnegative
integer $l$ such that $\sigma_d(X)=\sigma_d(Y)$
for $0\le d<l$ and $\sigma_l(X)<\sigma_l(Y)$.
If $X<Y$ or $X=Y$, then we shall write $X\le Y$.
Clearly this defines a total order on
the set of subsets of $\bZ_v$.

For $\pX=(X_1,\ldots,X_m)\in\pF_m$ and
any nonnegative integer $d$ we set
\[ \sigma_d(\pX)=\sigma_d(X_1)+\cdots+\sigma_d(X_m). \]

For $\pX,\pY\in\pF_m$, we shall write $\pX<\pY$ if there exists
a nonnegative integer $l$ such that $\sigma_d(\pX)=\sigma_d(\pY)$
for $1\le d<l$ and $\sigma_l(\pX)<\sigma_l(\pY)$.
If $\pX<\pY$ or $\pX=\pY$, then we shall write $\pX\le \pY$.
The binary relation $\le$ on $\pF_m$ is reflexive and transitive.
Moreover, if $\pX,\pY\in\pF_m$ then at least one of the
inequalities $\pX\le\pY$ and $\pY\le\pX$ holds.

We can now define our normal form.

\begin{definition} \label{normalna}
Let $\pX\in\pF_m$ and let $\pE\subseteq\pF_m$ be its equivalence class.
If there exists a $\pY\in\pE$ such that
\begin{itemize}
\item[(i)] $\pY\le\pZ$ for all $\pZ\in\pE$,
\item[(ii)] $\pY$ is unique up to a permutation,
\end{itemize}
then we say that $\pY$ is the \emph{normal form} of $\pX$.
We shall denote this normal form by $\nu(\pX)$.
\end{definition}

Let us point out some peculiar features of this normal form.

First of all, due to the condition (ii), it is not clear whether this
normal form always exists. 
It does in all cases that we have encountered so far.
It would be of interest to prove that (ii) is a consequence of (i).

Our normal form exhibits a very strong bias towards the small integers.
In particular each nonempty subset in the normal form contains 0.
All difference families in this paper will be given in this normal
form. To illustrate the normal form in the simplest case,
when $m=1$, we shall give in Table 1 the normal forms of all
cyclic $(v,k,\la)$ difference sets with $v\le50$.

\begin{remark} \label{dif-sets}
One can use cyclic difference sets to construct
difference families. For instance, $\pX=(\{0,1,4,6\},\{0,1,4,6\})$
is a $(13;4,4;2)$ difference family. By replacing the second
base block with the equivalent difference set
$\{0,2,8,12\}=2\cdot\{0,1,4,6\}$, we obtain the difference
family $\pY=(\{0,1,4,6\},\{0,2,8,12\})$ having the same parameters
$(13;4,4;2)$. Its normal form is
$\nu(\pY)=(\{0,1,4,6\},\{0,2,3,7\})$.
It is easy to check that these two difference families are
equivalent neither according to the definition \ref{equiv}
(e.g., they have different normal forms)
nor the one adopted in \cite{KK1}.
Hence, the claim made there that there is only one equivalence
class of $(13;4,4;2)$ difference families is in error.
This example shows that both definitions are inadequate when
dealing with the difference families having additional symmetry
properties. One can refine the notion of equivalence
and normal form to handle also such cases but we shall
not attempt to do it in this paper.
\end{remark}

\begin{tabular}{ll}
&\\
\multicolumn{2}{c}
{\bf  Table 1: Cyclic $(v,k,\la)$ difference sets with $v\le50$} \\
&\\
$(v,k,\la)$ & Difference set in the normal form \\
&\\ \hline
&\\
(7,3,1) & \{0,1,3\} \\
(11,5,2) & \{0,1,2,4,7\} \\
(13,4,1) & \{0,1,4,6\} \\
(15,7,3) & \{0,1,2,4,5,8,10\} \\
(19,9,4) & \{0,1,2,3,5,7,12,13,16\} \\
(21,5,1) & \{0,3,4,9,11\} \\
(23,11,5) & \{0,1,2,3,5,7,8,11,12,15,17\} \\
(31,6,1) & \{0,1,3,8,12,18\} \\
(31,15,7) & \{0,1,3,4,7,8,9,10,12,15,17,19,20,21,25\} \\
(35,17,8) & \{0,1,3,4,5,6,8,10,13,14,16,17,19,23,24,25,31\} \\
(37,9,2) & \{0,1,4,6,10,15,17,18,25\} \\
(40,13,4) & \{0,1,2,4,5,8,13,14,17,19,24,26,34\} \\
(43,21,10) & \{0,1,2,4,8,9,10,11,12,14,15,16,19,21,24,27,28,30,32, \\
 & \quad 33,37\} \\
(47,23,11) & \{0,1,2,3,5,6,7,8,11,13,15,16,17,20,23,24,26,27,31,33, \\
 & \quad 35,36,41\}
\end{tabular}

By using this normal form, we found that some families or
their equivalents have been listed more than once in the
same paper or by another author. For instance, in the paper
\cite{CK} one finds a list of four SDSs $(B_i,D_i)$
with parameters (31;15,10;10).
The last three of them are all equivalent:
\[ B_3=3B_2,\quad D_3=3D_2,\quad B_4=7B_3,\quad D_4=7D_3.
\quad\pmod{31} \]
The first SDS, $(B_1,D_1)$, has been rediscovered in \cite{MG}.

There are 10 SDSs in \cite{CK} with parameters (43;21,15;15).
All of them are indeed non-equivalent. The first one of them
again has been rediscovered in \cite{MG}.

We also found a few errors. For instance, in the paper \cite{CCS}
in the cases $(25;4,12;6)$ and $(27;3,5;1)$ the base blocks $D$
have wrong size, and in the case $(49;13,21;12)$ both blocks
$C$ and $D$ have wrong size. We have constructed the SDSs
in all three cases (see our tables below).

The often quoted table by Takeuchi \cite{KT} also has an
error: The item No. 31 is supposed to give an SDS with
parameters (17;8,8;7) but the first base block given there
is of size 7.

\section{Difference families with two base blocks}

From now on we shall consider only difference families with
exactly two base blocks (not necessarily distinct) and we shall
write their parameter sets as $(v;r,s;\la)$.
These parameters must satisfy the condition
\begin{equation} \label{jedn}
r(r-1)+s(s-1)=\la(v-1).
\end{equation}
We refer to parameter sets satisfying this necessary condition
as \emph{feasible}. There exist feasible parameter sets for which
there are no difference families. The first such example is
$(18;9,6;6)$, due to Young \cite{CY2}.
This also follows from a non-existence criterion of Arasu and
Xiang \cite{AX}.

There are a few infinite families of SDSs. One of them is due to
Szekeres, see e.g. \cite[p. 152]{AP},
another one to Koukouvinos, Kounias and Seberry \cite{KKS},
and several other families have been constructed by Wilson, Seberry, 
and Ding \cite{RW,JS,CD} using cyclotomy.

Clearly we may assume that $r\ge s$. Omitting the cases $s=0,1$
as trivial, we shall consider only the parameter sets satisfying
the inequalities
\begin{equation} \label{nejed}
v\ge2r\ge2s\ge4.
\end{equation}

For $v\le50$ there are exactly 227 feasible parameter sets
$(v;r,s;\la)$ satisfying these inequalities.
The difference families in this range with odd $v$ have been
constructed in many cases in the paper \cite{CCS}.

This work has been continued, including the even $v$'s, in
the papers \cite{DZ1,MDV}. The existence question has been
resolved for all feasible parameter sets with $v\le40$,
but in the range $40<v\le50$ there remained 23 undecided
cases. We shall settle eight of these cases by
constructing the required difference families.

In the cases that we consider, the order is given by $n=r+s-\la$.
We mention that it plays an important role in the criterion of
Arasu and Xiang. We shall now introduce three special types of
SDSs which have been studied extensively.
Before we do that, let us recall some basic definitions.

Let $\pA$ be a binary sequence of length $v$, i.e.,
$\pA=(a_0,a_1,\ldots,a_{v-1})$ where each $a_i$ is $\pm1$.
The \emph{periodic autocorrelation function} (PACF)
$\tilde{\phi}$ of $\pA$ is defined by
\[ \tilde{\phi}(i)=\sum_{j=0}^{v-1} a_j a_{i+j},
\quad 0\le i<v, \]
where $i+j$ should be reduced modulo $v$.
The \emph{non-periodic autocorrelation function} (NACF)
$\phi$ of $\pA$ is defined by
\[ \phi(i)=\sum_{j=0}^{v-1-i} a_j a_{i+j},
\quad 0\le i<v. \]

The first type is defined by the condition $r=s$, i.e., the
two base blocks are of the same size. Such SDSs give
BIBD's on $\bZ_v$ by an old result of Bose \cite{RB}.

The second type is defined by the relation $v=2n+1$.
They are known as D-optimal SDSs because they can be used
to construct D-optimal designs (of circular type) of size $2v$.
The D-optimal parameter sets in our range are listed in Table 2.
The columns with heading ``\#'' give the number of equivalence classes
of SDSs. If this number is not known, the question mark is entered.

The equivalence classes of the D-optimal SDSs have been
enumerated in \cite{CY1} for $v\le19$, and later this has been
extended in \cite{KK1,KK2} to all $v\le27$ and $v=33,45$.

\begin{tabular}{lr|lr|lr}
\\
\multicolumn{6}{c}
{\bf  Table 2: D-optimal parameters $(v;r,s;\la)$} \\
\\
 & \# &  & \# &  & \# \\ \hline
(9;3,2;1) & 1 & \quad(13;4,4;2) & 2 & \quad(13;6,3;3) & 2 \\
(15;6,4;3) & 3 & \quad(19;7,6;4) & 8 & \quad(21;10,6;6) & 31 \\
(23;10,7;6) & 17 & \quad(25;9,9;6) & 39 & \quad(27;11,9;7) & 48 \\
(31;15,10;10) & ? & \quad(33;13,12;9) & 509 & \quad(33;15,11;10) & 516 \\
(37;16,13;11) & ? & \quad(41;16,16;12) & ? & \quad(43;18,16;13) & ? \\
(43;21,15;15) & ? & \quad(45;21,16;15) & 1358 & \quad(49;22,18;16) & ? \\
\multicolumn{6}{c}{} \\
\end{tabular}

The third type of SDSs is defined by the relation $v=2n$.
They are essentially the same objects as the pairs $\pA=(A_1,A_2)$
of binary sequences of length $v$ for which the sum of
the PACF of $A_1$ and $A_2$
is zero (more precisely, a $\delta$-function).
We shall refer to this sum as the
\emph{periodic autocorrelation function} of $\pA$.
Due to their importance in engineering applications,
such binary sequences have been studied for long time
(see e.g. \cite{FSX} and its references).

If one replaces above the word ``periodic'' by ``aperiodic'' one
obtains the definition of Golay pairs.
We refer to the lengths $v$ of Golay pairs as \emph{Golay numbers}.
The equivalence classes of
Golay pairs have been described in our paper \cite{DZ3}
and the representatives of these classes given for Golay numbers
2,4,8,10,16,20,26,32 and 40.
For further results in this direction see \cite{BF}.
Each Golay pair $(A_1,A_2)$ of length $v$ gives a cyclic
$(v;r,s;\la)$ difference family $(X_1,X_2)$ 
by taking $X_i$ to be the set of positions of terms $-1$
in the sequence $A_i$. Here is an example.

\begin{example} \label{G-par}
Let us start with the Golay pair of length 10:
\begin{eqnarray*}
A_1 &=& +,+,-,+,-,+,-,-,+,+; \\
A_2 &=& +,+,-,+,+,+,+,+,-,-.
\end{eqnarray*}
We label the positions by integers 0 to 9 from left
to right. The positions of the negative signs (which stand
for $-1$) are recorded in the two base blocks
$X_1 = \{2,4,6,7\}$ and $X_2 = \{2,8,9\}$.
Then $(X_1,X_2)$ is a difference family in $\bZ_{10}$.
Its normal form $(\{0,1,3,5\},\{0,1,4\})$ appears in Table 3.
\end{example}

\begin{remark}
There is a natural definition of equivalence
for Golay pairs, see e.g. \cite{BF,DZ3}.
While equivalent Golay pairs always give equivalent SDSs,
two non-equivalent Golay pairs may also give rise to
equivalent SDSs.
For instance, there are 5 equivalence classes of Golay
pairs of length 8  but they give only 2 equivalence
classes of SDSs.
\end{remark}

We shall split the table of difference families into four
parts, Tables 3-6 covering the ranges $v\le20$,
$20<v\le30$, $30<v\le40$ and $40<v\le50$, respectively.
The parameter sets are arranged in the lexicographic
order by the parameters $v$, $\la$ and $n$.
The asterisk in the last column of the tables means that
the family has been constructed by using our genetic type
program and that we were not able to locate an
equivalent family in the literature.

\section{The range $v\le20$}

In Table 3 we list the 30 feasible parameter sets $(v;r,s;\la)$
with $v\le20$ satisfying (\ref{nejed}) and give
examples of cyclic difference families for them.
In four of these cases the difference families do not exist.
This is indicated
in Table 3 (and subsequent tables) by the word ``None''.

In each case we provide in the last column
a reference where the family has
been constructed (or just listed) or where it was established
that the families do not exist. No effort has been made to
assign priorities for these results.
The symbol ``Sz'' indicates that the difference family is
a classical Szekeres family, and ``DS'' means that the two
blocks are in fact difference sets.
Whenever possible, we have tried to avoid using the
``DS'' type examples.

The symbol ``GP'' indicates that the difference family is
obtained from a Golay pair.
In all four cases $v=8,10,16,20$ we have included in the table all
equivalence classes of SDSs arising from Golay pairs of length $v$.

\newpage
\begin{tabular}{lrlll}
\multicolumn{5}{c}
{\bf Table 3: $(v;r,s;\la)$ difference families with $v\le20$} \\
&&&&\\
\multicolumn{1}{l}{$(v;r,s;\la)$} &
\multicolumn{1}{l}{$n$} &
\multicolumn{2}{c}{Base blocks} &
\multicolumn{1}{l}{Ref.} \\
&&&&\\
\hline
&&&&\\
(5;2,2;1) & 3 & \{0,1\} & \{0,2\} & \cite{JS} \\
(7;3,3;2) & 4 & \{0,1,3\} & \{0,1,3\} & DS \\
(8;4,2;2) & 4 & \{0,1,2,4\} & \{0,3\} & GP \\
 && \{0,1,3,4\} & \{0,2\} & GP \\
(9;3,2;1) & 4 & \{0,1,4\} & \{0,2\} & \cite{HE} \\
(9;4,4;3) & 5 & \{0,1,3,4\} & \{0,1,3,5\} & Sz \\
(10;4,3;2) & 5 & \{0,1,3,5\} & \{0,1,4\} & GP \\
(11;5,5;4) & 6 & \{0,1,2,4,6\} & \{0,1,2,5,8\} & Sz \\
(12;5,2;2) & 5 & \{0,1,2,5,8\} & \{0,2\} & \cite{DZ1} \\
(13;3,3;1) & 5 & \{0,1,4\} & \{0,2,7\} & \cite{KT} \\
(13;4,4;2) & 6 & \{0,1,4,6\} & \{0,1,4,6\} & DS \\
 && \{0,1,4,6\} & \{0,2,3,7\} & DS \\
(13;6,3;3) & 6 & \{0,1,2,4,7,9\} & \{0,1,4\} & \cite{JS} \\
 && \{0,1,3,5,7,8\} & \{0,1,4\} & \cite{HE} \\
(13;6,6;5) & 7 & \{0,1,2,3,6,9\} & \{0,1,2,4,6,9\} & \cite{KT} \\
 && \{0,1,2,3,6,10\} & \{0,1,3,5,7,8\} & \cite{JS} \\
 && \{0,1,2,4,5,8\} & \{0,1,2,4,7,9\} & \cite{JS} \\
 && \{0,1,2,4,5,8\} & \{0,1,3,5,7,8\} &$\ast$ \\
(14;5,3;2) & 6 & None & & \cite{DZ1} \\
(15;4,2;1) & 5 & \{0,1,4,9\} & \{0,2\} & \cite{CCS} \\
(15;6,4;3) & 7 & \{0,1,2,4,6,9\} & \{0,1,4,9\} & \cite{HE} \\
 && \{0,1,3,4,8,10\} & \{0,1,4,6\} & \cite{CY1} \\
 && \{0,1,3,5,7,8\} & \{0,1,4,10\} & \cite{CY1} \\
(15;7,7;6) & 8 & \{0,1,2,3,4,8,11\} & \{0,1,3,5,6,9,11\} & Sz \\
 && \{0,1,2,4,5,7,9\} & \{0,1,2,5,7,8,11\} &$\ast$ \\
 && \{0,1,2,4,5,9,11\} & \{0,1,3,4,6,8,9\} &$\ast$ \\
(16;6,6;4) & 8 & \{0,1,2,3,6,10\} & \{0,1,3,6,8,12\} & GP \\
 && \{0,1,2,4,5,10\} & \{0,1,4,7,9,11\} & GP \\
 && \{0,1,2,4,6,9\} & \{0,1,2,6,9,12\} & GP \\
 && \{0,1,2,4,6,9\} & \{0,1,5,7,8,11\} & GP \\
 && \{0,1,2,5,9,11\} & \{0,1,3,5,6,9\} &$\ast$ \\
 && \{0,1,3,4,7,9\} & \{0,1,5,6,8,10\} & GP \\
 && \{0,1,3,5,6,9\} & \{0,1,3,5,9,10\} & GP \\
 && \{0,1,3,5,7,8\} & \{0,1,4,6,9,10\} & GP \\
(17;6,2;2) & 6 & None & & \cite{DZ1} \\
\end{tabular}

\newpage
\begin{tabular}{lrlll}
\multicolumn{5}{c}{\bf Table 3 (continued)} \\
&&&&\\
\multicolumn{1}{l}{$(v;r,s;\la)$} &
\multicolumn{1}{l}{$n$} &
\multicolumn{2}{c}{Base blocks} &
\multicolumn{1}{l}{Ref.} \\
&&&&\\
\hline
&&&&\\
(17;5,4;2) & 7 & \{0,1,4,6,10\} & \{0,1,3,8\} &$\ast$ \\
 && \{0,1,4,7,9\} & \{0,1,5,7\} & \cite{CCS} \\
(17;7,3;3) & 7 & None & & \cite{DZ1} \\
(17;8,8;7) & 9 & \{0,1,2,4,5,7,10,11\} & \{0,1,2,4,6,8,9,14\} &$\ast$ \\
 && \{0,2,3,4,7,8,9,11\} & \{0,2,3,5,6,8,12,13\} & \cite{JS} \\
(18;8,4;4) & 8 & \{0,1,2,3,5,8,9,13\} & \{0,2,6,9\} & \cite{DZ1} \\
 && \{0,1,2,4,6,7,10,11\} & \{0,3,5,11\} &$\ast$ \\
(18;9,6;6) & 9 & None & & \cite{AX,CY2} \\
(19;4,3;1) & 6 & \{0,3,5,9\} & \{0,1,8\} & \cite{CCS} \\
(19;6,3;2) & 7 & \{0,1,2,6,10,13\} & \{0,2,5\} &$\ast$ \\
 && \{0,1,2,7,8,11\} & \{0,2,5\} & \cite{CCS} \\
(19;7,4;3) & 8 & \{0,1,2,4,7,9,13\} & \{0,1,4,9\} &$\ast$ \\
 && \{0,1,3,4,8,10,14\} & \{0,1,3,8\} & \cite{CCS} \\
(19;7,6;4) & 9 & \{0,1,2,3,7,11,14\} & \{0,2,5,6,9,11\} & \cite{CY1} \\
 && \{0,1,2,4,5,10,13\} & \{0,1,4,6,8,13\} & \cite{CY1} \\
 && \{0,1,2,5,7,11,12\} & \{0,2,4,5,8,11\} & \cite{CY1} \\
 && \{0,1,2,5,7,11,14\} & \{0,2,3,4,8,11\} & \cite{HE} \\
 && \{0,1,3,4,7,12,14\} & \{0,1,2,5,7,11\} & \cite{HE} \\
 && \{0,1,3,4,8,10,14\} & \{0,1,2,4,7,12\} & \cite{CY1} \\
 && \{0,1,3,4,8,10,14\} & \{0,2,3,7,8,10\} & \cite{CY1} \\
 && \{0,2,3,4,6,9,14\} & \{0,1,2,6,10,13\} & \cite{CY1} \\
(19;9,9;8) & 10 & \{0,1,2,3,4,8,10,13,15\} & \{0,1,2,4,5,6,9,12,15\}
 & \cite{DZ1} \\
 && \{0,1,2,4,6,8,9,11,12\} & \{0,1,2,5,6,7,9,12,15\} &$\ast$ \\
(20;8,5;4) & 9 & \{0,1,2,4,5,8,11,13\} & \{0,2,6,7,12\} &$\ast$ \\
 && \{0,1,2,5,6,8,11,13\} & \{0,1,3,7,11\} & \cite{DZ1} \\
(20;9,7;6) & 10 & \{0,1,2,3,4,6,9,10,14\} & \{0,1,4,8,10,13,15\} & GP \\
 && \{0,1,2,3,5,6,9,13,15\} & \{0,2,3,5,9,10,14\} & GP \\
 && \{0,1,2,3,5,7,9,12,13\} & \{0,1,3,6,7,10,15\} & GP \\
 && \{0,1,2,3,5,7,10,11,14\} & \{0,1,3,5,8,9,15\} & GP \\
 && \{0,1,2,3,6,7,9,11,15\} & \{0,2,3,5,9,10,13\} &$\ast$ \\
 && \{0,1,2,3,6,7,10,12,14\} & \{0,1,4,5,7,10,12\} & GP \\
 && \{0,1,2,3,6,7,11,13,15\} & \{0,2,3,5,8,9,12\} & GP \\
 && \{0,1,2,3,6,8,10,11,15\} & \{0,1,3,4,7,9,13\} & GP \\
 && \{0,1,2,3,6,9,11,13,17\} & \{0,1,2,5,7,8,12\} & GP \\
 && \{0,1,2,4,6,7,9,10,14\} & \{0,1,2,5,9,11,14\} & GP \\
 && \{0,1,2,4,7,9,10,11,15\} & \{0,1,3,5,6,9,13\} &$\ast$ \\
\end{tabular}

\newpage
\begin{tabular}{lrlll}
\multicolumn{5}{c}{\bf Table 3 (continued)} \\
&&&&\\
\multicolumn{1}{l}{$(v;r,s;\la)$} &
\multicolumn{1}{l}{$n$} &
\multicolumn{2}{c}{Base blocks} &
\multicolumn{1}{l}{Ref.} \\
&&&&\\
\hline
&&&&\\
(20;9,7;6) & 10 & \{0,1,2,5,6,7,10,12,14\} & \{0,1,3,4,7,10,12\} & GP \\
 && \{0,1,2,5,7,8,10,12,16\} & \{0,1,2,3,6,9,13\} & GP \\
 && \{0,1,3,4,5,7,9,13,14\} & \{0,1,2,5,7,10,13\} & GP \\
 && \{0,1,3,4,6,8,9,13,14\} & \{0,1,2,4,6,10,13\} & GP \\
 && \{0,1,3,4,6,9,11,13,17\} & \{0,1,2,3,7,8,12\} & GP \\
 && \{0,1,3,4,7,8,9,12,14\} & \{0,1,2,4,6,11,14\} & \cite{DZ1} \\
 && \{0,2,3,6,7,9,10,12,14\} & \{0,1,2,5,6,11,13\} & GP \\
&&&&\\
\end{tabular}

The difference families in the five D-optimal cases in this
range have been enumerated, up to equivalence, in \cite{CY1}.
As there are not many of them, we have included all
of them in the table.

As a curiosity, let us mention that Ehlich gives in his paper
\cite{HE} three examples of D-optimal designs of circular type and 
size $2v=38$. Among the corresponding three SDSs, all with
parameters (19;7,6;4), the first two are equivalent.
Similarly, his two examples in the case $2v=18$ are equivalent.
(In the case $2v=26$, his two SDSs are indeed non-equivalent.)

\section{The range $20<v\le30$}

There are 45 feasible parameter sets in this range.
Only four of them are D-optimal (marked by the symbol ``DO'').
The representatives of the equivalence classes of these
D-optimal SDSs have been computed in \cite{KK1}.
The number of classes is 31, 17, 39 and 48 when
$v$ is 21, 23, 25 and 27, respectively.
We do not list them and refer the reader to that paper.

There is only one Golay number, $v=26$, in this range
and there is up to equivalence a unique Golay pair
of that length The corresponding difference family with
parameters $(26;11,10;8)$ is included in the table
(and marked with ``GP'').

\newpage

\begin{tabular}{ll}
\multicolumn{2}{c}
{\bf Table 4: $(v;r,s;\la)$ difference families with $20<v\le30$} \\
 & \\
\hline
 & \\
$(21;5,5;2),\quad n=8$\quad None & \cite{DZ1} \\
 & \\
$(21;6,6;3),\quad n=9$ & \\
\quad\{0,1,3,7,10,15\},\quad\{0,2,3,4,8,13\}; & \cite{KT} \\
\quad\{0,1,3,8,10,14\},\quad\{0,1,4,6,9,10\}; &$\ast$ \\
 & \\
$(21;10,6;6),\quad n=10$ & DO\\
 & \\
$(21;10,10;9),\quad n=11$ & \\
\quad\{0,1,2,3,4,7,8,10,13,14\},\quad\{0,2,4,5,7,8,10,12,16,17\}; & Sz \\
\quad\{0,1,2,3,6,7,9,11,12,14\},\quad\{0,1,3,4,5,8,9,11,15,17\}; &$\ast$ \\
\quad\{0,1,2,4,5,6,7,10,13,15\},\quad\{0,1,2,4,7,8,11,12,14,16\};
 & \cite{HCD} \\
 & \\
$(22;6,4;2),\quad n=8$ & \\
\quad\{0,1,4,6,10,15\},\quad\{0,2,3,10\}; & \cite{DZ1} \\
 & \\
$(22;9,4;4),\quad n=9$ & \\
\quad\{0,1,2,4,6,10,11,14,17\},\quad\{0,1,3,8\}; &$\ast$ \\
\quad\{0,1,2,5,7,8,10,12,16\},\quad\{0,1,4,13\}; &$\ast$ \\
\quad\{0,1,2,5,7,9,11,12,15\},\quad\{0,1,6,9\}; & \cite{DZ1} \\
 & \\
$(22;7,7;4),\quad n=10$ & \\
\quad\{0,1,2,3,7,9,12\},\quad\{0,3,4,8,10,14,17\}; & \cite{HCD} \\
\quad\{0,1,2,3,7,9,13\},\quad\{0,2,5,6,9,14,17\}; &$\ast$ \\
\quad\{0,1,2,4,8,10,13\},\quad\{0,1,5,6,8,11,15\};&$\ast$ \\
\quad\{0,1,3,4,8,11,17\},\quad\{0,1,3,5,7,12,13\}; & \cite{DZ1} \\
 & \\
$(23;5,2;1),\quad n=6$ & \\
\quad\{0,1,3,8,14\},\quad\{0,4\}; & \cite{CCS}\\
 & \\
$(23;7,2;2),\quad n=7$\quad None & \cite{DZ1} \\
 & \\
$(23;10,5;5),\quad n=10$ & \\
\quad\{0,1,2,3,5,7,8,12,13,16\},\quad\{0,2,6,9,15\};&$\ast$ \\
\quad\{0,1,2,3,5,7,8,12,14,18\},\quad\{0,1,2,8,11\};&$\ast$ \\
\quad\{0,1,2,4,6,10,11,12,15,18\},\quad\{0,2,3,7,10\}; & \cite{DZ1} \\
 & \\
$(23;10,7;6),\quad n=11$ & DO\\
\end{tabular}

\begin{tabular}{ll}
\multicolumn{2}{c} {\bf Table 4 (continued)} \\
 & \\
\hline
 & \\
$(23;11,11;10),\quad n=12$ & \\
\quad\{0,1,2,3,4,6,9,10,14,15,17\},\quad\{0,1,3,5,6,7,9,11,12,16,19\};& Sz \\
\quad\{0,1,2,3,6,9,10,11,13,15,18\},\quad\{0,1,2,4,6,7,8,11,12,14,20\};&$\ast$ \\
 & \\
$(24;10,2;4),\quad n=8$ & \\
\quad\{0,1,2,3,5,8,9,13,15,18\},\quad\{0,4\};&$\ast$ \\
\quad\{0,1,2,5,7,8,11,12,14,16\},\quad\{0,8\};&$\ast$ \\
\quad\{0,1,4,5,6,7,10,13,15,17\},\quad\{0,8\};& \cite{DZ1} \\
 & \\
$(24;9,5;4),\quad n=10$ & \\
\quad\{0,1,2,3,5,8,10,14,18\},\quad\{0,1,6,10,13\};&$\ast$ \\
\quad\{0,1,2,4,5,7,11,12,16\},\quad\{0,2,8,11,18\};&$\ast$ \\
\quad\{0,1,2,4,5,10,12,16,19\},\quad\{0,2,5,6,13\};&$\ast$ \\
\quad\{0,1,2,4,5,11,12,16,19\},\quad\{0,2,6,8,11\};&$\ast$ \\
\quad\{0,1,2,5,6,9,11,16,18\},\quad\{0,2,3,6,14\};&$\ast$ \\
\quad\{0,1,3,4,8,10,12,13,18\},\quad\{0,2,5,6,13\};&$\ast$ \\
\quad\{0,1,3,5,9,10,11,14,17\},\quad\{0,1,3,7,12\};& \cite{DZ1} \\
 & \\
$(24;12,3;6),\quad n=9$ & \\
\quad\{0,1,2,5,6,8,9,11,13,14,15,18\},\quad\{0,2,10\};& \cite{DZ1} \\
 & \\
$(25;4,4;1),\quad n=7$\quad None & \cite{DZ1} \\
$(25;7,3;2),\quad n=8$\quad None & \cite{DZ1} \\
 & \\
$(25;7,6;3),\quad n=10$ & \\
\quad\{0,1,2,5,8,12,14\},\quad\{0,1,5,8,10,16\};& \cite{CCS} \\
\quad\{0,1,3,8,9,12,15\},\quad\{0,2,3,7,8,12\};&$\ast$ \\
 & \\
$(25;10,3;4),\quad n=9$\quad None & \cite{DZ1} \\
 & \\
$(25;10,6;5),\quad n=11$ & \\
\quad\{0,1,3,4,5,8,11,15,17,20\},\quad\{0,1,2,6,8,15\};&$\ast$ \\
\quad\{0,1,3,4,5,9,10,12,16,18\},\quad\{0,1,4,6,11,14\};& \cite{CCS} \\
 & \\
$(25;12,4;6),\quad n=10$ & \\
\quad\{0,1,2,4,5,6,8,10,13,15,18,19\},\quad\{0,1,7,10\};&$\ast$ \\
 & \\
$(25;9,9;6),\quad n=12$ & DO\\
\end{tabular}

\begin{tabular}{ll}
\multicolumn{2}{c} {\bf Table 4 (continued)} \\
 & \\
\hline
 & \\
$(25;12,12;11),\quad n=13$ & \\
\quad\{0,1,2,3,5,7,11,12,14,17,19,20\}, & \\
\quad\{0,1,2,4,5,6,8,10,11,14,15,18\};& \cite{DZ1} \\
\quad\{0,1,2,3,7,8,11,12,13,15,17,20\}, & \\
\quad\{0,1,2,4,5,7,8,10,11,14,16,18\};&$\ast$ \\
 & \\
$(26;6,5;2),\quad n=9$ & \\
\quad\{0,1,2,8,13,17\},\quad\{0,2,5,8,12\};& \cite{DZ1} \\
\quad\{0,2,3,8,11,15\},\quad\{0,1,5,7,17\};&$\ast$ \\
 & \\
$(26;11,10;8),\quad n=13$ & \\
\quad\{0,1,2,3,4,7,9,12,14,16,20\},\quad\{0,1,4,5,9,10,11,13,16,19\}; &
\cite{DZ1} \\
\quad\{0,1,2,4,5,8,10,14,16,19,21\},\quad\{0,1,2,4,5,8,11,12,13,18\}; & GP \\
\quad\{0,2,3,4,5,7,10,11,14,18,20\},\quad\{0,1,2,5,6,7,11,14,17,19\};&$\ast$ \\
 & \\
$(27;5,3;1),\quad n=7$ & \\
\quad\{0,2,5,11,15\},\quad\{0,1,8\};&$\ast$ \\
 & \\
$(27;9,3;3),\quad n=9$\quad None & \cite{AX,MDV} \\
 & \\
$(27;11,5;5),\quad n=11$ & \\
\quad\{0,1,2,4,5,8,10,11,16,18,23\},\quad\{0,1,4,11,13\};&$\ast$ \\
\quad\{0,1,2,5,6,7,10,13,16,18,20\},\quad\{0,3,4,10,12\};
 & \cite{CCS} \\
 & \\
$(27;11,9;7),\quad n=13$ & DO\\
 & \\
$(27;13,13;12),\quad n=14$ & \\
\quad\{0,1,2,3,4,6,7,8,13,14,17,19,22\}, & \\
\quad\{0,1,3,4,5,8,10,11,12,14,18,20,23\};&$\ast$ \\
\quad\{0,1,2,3,5,6,8,9,12,16,17,19,21\}, & \\
\quad\{0,1,3,4,5,7,9,10,13,14,15,20,22\};& \cite{DZ1} \\
 & \\
$(28;7,4;2),\quad n=9$ & \\
\quad\{0,2,3,5,9,13,19\},\quad\{0,1,8,13\};& \cite{DZ1} \\
 & \\
$(28;13,3;6),\quad n=10$ & \\
\quad\{0,1,2,3,6,7,8,10,12,15,17,18,21\},\quad\{0,4,12\};&$\ast$ \\
\quad\{0,1,3,4,6,7,11,12,13,16,18,20,21\},\quad\{0,2,6\};& \cite{DZ1} \\
\end{tabular}

\begin{tabular}{ll}
\multicolumn{2}{c} {\bf Table 4 (continued)} \\
 & \\
\hline
 & \\
$(28;12,6;6),\quad n=12$\quad None & \cite{MDV} \\
 & \\
$(28;10,9;6),\quad n=13$ & \\
\quad\{0,1,2,5,6,7,9,14,17,19\},\quad\{0,1,5,8,11,12,14,18,20\};
 & \cite{DZ1} \\
\quad\{0,2,3,4,5,9,13,15,16,23\},\quad\{0,2,4,5,9,10,13,16,22\};& $\ast$ \\
 & \\
$(29;9,4;3),\quad n=10$ & \\
\quad\{0,1,2,4,7,8,12,18,21\},\quad\{0,2,7,16\};& \cite{DZ1} \\
 & \\
$(29;7,7;3),\quad n=11$ & \\
\quad\{0,1,2,5,9,18,24\},\quad\{0,2,4,7,10,18,19\};& \cite{KT} \\
\quad\{0,1,4,6,11,12,20\},\quad\{0,1,4,7,12,14,16\};& $\ast$ \\
 & \\
$(29;11,2;4),\quad n=9$\quad None & \cite{DZ1} \\
 & \\
$(29;8,8;4),\quad n=12$ & \\
\quad\{0,1,2,4,8,10,15,20\},\quad\{0,1,4,7,8,13,16,18\};& $\ast$ \\
\quad\{0,1,5,6,7,10,14,23\},\quad\{0,2,4,7,10,18,19,21\};& \cite{KT}\\
 & \\
$(29;11,6;5),\quad n=12$ & \\
\quad\{0,1,2,3,4,8,9,13,17,20,23\},\quad\{0,2,4,7,12,18\};& \cite{DZ1} \\
\quad\{0,1,2,3,6,9,11,13,17,20,25\},\quad\{0,1,4,6,13,14\};&$\ast$ \\
 & \\
$(29;13,4;6),\quad n=11$ & \\
\quad\{0,1,2,3,6,7,10,11,13,15,16,18,22\},\quad\{0,2,8,19\};
 & \cite{DZ1} \\
 & \\
$(29;14,7;8),\quad n=13$ & \\
\quad\{0,1,2,5,6,7,9,11,13,14,15,18,21,24\},
\quad\{0,2,3,4,11,14,19\};&$\ast$ \\
\quad\{0,2,3,4,6,7,9,10,11,13,18,19,23,24\},
\quad\{0,2,4,7,10,18,19\};& \cite{CCS} \\
 & \\
$(29;14,14;13),\quad n=15$ & \\
\quad\{0,1,2,3,4,7,9,11,12,15,17,19,20,26\}, & \\ 
\quad\{0,1,2,4,5,6,8,9,11,14,18,22,23,24\};& \cite{JS} \\
\quad\{0,1,2,3,4,7,9,11,12,15,17,19,20,26\}, & \\ 
\quad\{0,2,3,4,6,7,9,10,11,13,18,19,23,24\};& \cite{CD} \\
\quad\{0,1,2,3,4,7,9,12,13,16,17,19,21,22\}, & \\ 
\quad\{0,1,2,4,5,7,8,10,11,12,16,18,23,25\};& Sz \\
\quad\{0,1,2,3,5,6,10,12,13,15,16,17,21,23\}, & \\ 
\quad\{0,1,2,4,5,7,9,10,11,13,16,17,21,24\};&$\ast$ \\
\end{tabular}

\begin{tabular}{ll}
\multicolumn{2}{c} {\bf Table 4 (continued)} \\
 & \\
\hline
 & \\
$(29;14,14;13),\quad n=15$ & \\
\quad\{0,1,2,4,5,6,9,11,13,14,15,17,21,24\}, & \\ 
\quad\{0,1,3,4,6,8,9,11,12,16,17,18,19,23\};&$\ast$ \\
\quad\{0,2,3,4,5,8,10,13,14,15,16,18,22,25\}, & \\ 
\quad\{0,2,3,4,6,7,9,10,11,13,18,19,23,24\};& \cite{JS} \\
 & \\
$(30;8,2;2),\quad n=8$\quad None & \cite{DZ1} \\
$(30;11,3;4),\quad n=10$\quad None & \cite{MDV} \\
 & \\
$(30;12,7;6),\quad n=13$ & \\
\quad\{0,1,2,4,5,7,8,9,14,17,19,23\},\quad\{0,2,6,11,12,19,22\};&$\ast$ \\
\quad\{0,1,2,4,5,8,10,14,17,19,23,24\},
\quad\{0,1,2,4,7,12,20\};& \cite{DZ1} \\
\quad\{0,1,3,6,7,8,10,11,13,18,22,24\},\quad\{0,1,4,5,13,15,21\};&$\ast$ \\
\end{tabular}

\section{The range $30<v\le40$}

In this range we have 63 feasible parameter sets, and only
two Golay numbers: 32 and 40. There are
330 equivalence classes of Golay pairs of length 32, and
220 of length 40. The difference families that arise from
these pairs have parameters $(32;16,12;12)$ and $(40;18,16;14)$
respectively (up to equivalence). We have included in Table 5
only one example of such difference families for each of
these parameter sets.

There is one more set of parameters, namely
$(34;16,13;12)$, which belongs to the special type defined
by the condition $v=2n$. This is an important case since
the periodic autocorrelation function of
the binary sequences $\pA=(A_1,A_2)$ constructed from the
SDSs $\pX=(X_1,X_2)$ with these parameters is zero.
We have constructed two non-equivalent such SDSs in our paper
\cite{DZ2}. In fact one of them appears in our earlier
paper \cite{DZ1} but at that time its significance was not
recognized. Now we have constructed an additional
non-equivalent SDS with the same parameters (see Table 5).

We also have in this range four D-optimal parameter sets.
The equivalence classes of SDSs for the two
sets with $v=33$ have been enumerated in \cite{KK1}.
We do not give any difference families in these two cases. 
For $v=31$ we have included the two non-equivalent D-optimal
SDSs from \cite{CK}, and for $v=37$ the two known
D-optimal SDSs \cite{JC,MG}.
We have contributed two new D-optimal designs for each
of the cases $v=31$ and $v=37$.

\newpage
\begin{tabular}{ll}
\multicolumn{2}{c}
{\bf Table 5: $(v;r,s;\la)$ difference families with $30<v\le40$} \\
 & \\
\hline
 & \\
$(31;6,6;2),\quad n=10$ & \\
\quad\{0,1,4,11,12,17\},\quad\{0,2,4,9,12,18\}; &$\ast$ \\
\quad\{0,1,4,11,13,17\},\quad\{0,3,5,10,11,19\}; & \cite{DZ1} \\
 & \\
$(31;10,6;4),\quad n=12$ & \\
\quad\{0,1,2,6,8,9,11,14,18,22\},\quad\{0,1,5,7,17,20\};&$\ast$ \\
\quad\{0,3,4,9,11,14,16,21,22,25\},\quad\{0,1,2,4,8,16\};& \cite{CCS} \\
 & \\
$(31;10,10;6),\quad n=14$ & \\
\quad\{0,1,2,3,5,7,11,15,16,23\},\quad\{0,1,4,6,7,11,14,20,23,25\};&
\cite{CCS} \\
\quad\{0,1,2,5,6,9,12,15,17,19\},\quad\{0,1,2,5,7,11,13,16,23,24\};&$\ast$ \\
 & \\
$(31;15,6;8),\quad n=13$ & \\
\quad\{0,1,2,3,4,5,9,10,12,15,19,21,23,26,27\},\quad\{0,2,3,6,13,18\};&
\cite{DZ1} \\
\quad\{0,1,2,3,5,7,9,10,12,13,16,18,21,22,26\},\quad\{0,1,3,7,14,15\};&$\ast$ \\
 & \\
$(31;15,10;10),\quad n=15$ & \\
\quad\{0,1,2,3,5,6,7,11,13,15,16,18,23,24,27\},& \\
\quad\{0,2,3,5,6,8,12,19,20,27\};& \cite{CK} \\
\quad\{0,1,2,3,5,7,8,9,11,14,17,18,19,22,27\},& \\
\quad\{0,1,2,6,8,9,12,16,19,21\};&$\ast$ \\
\quad\{0,1,2,4,5,6,7,9,11,14,15,17,21,23,24\},& \\
\quad\{0,1,3,5,6,11,12,16,19,23\};&$\ast$ \\
\quad\{0,1,2,5,6,8,9,11,12,14,18,21,23,25,26\},& \\
\quad\{0,1,2,3,5,7,11,15,16,23\};& \cite{CK} \\
 & \\
$(31;15,15;14),\quad n=16$ & \\
\quad\{0,1,2,3,4,5,7,12,13,14,17,20,22,26,28\},& \\
\quad\{0,1,2,3,6,7,9,10,12,14,16,17,20,21,25\};& \cite{DZ1} \\
\quad\{0,1,2,3,4,6,7,8,12,13,15,16,18,22,26\},& \\
\quad\{0,1,2,3,5,7,9,12,14,15,17,18,22,25,26\};&$\ast$ \\
 & \\
$(32;8,3;2),\quad n=9$\quad None & \cite{MDV} \\
 & \\
$(32;7,5;2),\quad n=10$ & \\
\quad\{0,2,3,5,10,14,18\},\quad\{0,1,7,13,22\};& \cite{DZ1} \\
 & \\
$(32;13,6;6),\quad n=13$ & \\
\quad\{0,1,2,3,4,8,9,11,13,15,18,19,23\},\quad\{0,3,6,12,19,24\};&$\ast$ \\
\quad\{0,1,3,4,5,9,10,12,14,17,20,22,26\},\quad\{0,1,3,7,14,15\};&
\cite{DZ1} \\
\end{tabular}

\begin{tabular}{ll}
\multicolumn{2}{c} {\bf Table 5 (continued)} \\
 & \\
\hline
 & \\
$(32;16,12;12),\quad n=16$ & \\
\quad\{0,1,2,3,6,7,10,11,12,14,16,18,20,23,25,26\},& \\
\quad\{0,1,2,3,6,8,11,13,14,17,18,21\};& GP \\
\quad\{0,1,2,4,5,6,7,10,12,13,15,17,18,22,24,26\},& \\
\quad\{0,1,3,4,5,9,12,13,16,18,19,26\};&$\ast$ \\
\quad\{0,1,2,4,6,7,8,10,12,13,15,17,18,21,22,26\},& \\
\quad\{0,1,2,4,5,9,11,12,18,19,21,24\};& \cite{DZ1} \\
 & \\
$(33;6,2;1),\quad n=7$\quad None & \cite{DZ1} \\
$(33;5,4;1),\quad n=8$\quad None & \cite{DZ1} \\
 & \\
$(33;10,3;3),\quad n=10$ & \\
\quad\{0,1,2,5,7,8,12,16,22,25\},\quad\{0,2,14\};& \cite{CCS} \\
 & \\
$(33;9,8;4),\quad n=13$ & \\
\quad\{0,1,2,7,10,11,14,17,22\},\quad\{0,2,4,6,11,14,19,20\};&$\ast$ \\
\quad\{0,2,5,7,8,14,17,18,22\},\quad\{0,1,2,6,9,13,15,23\};& \cite{MDV} \\
 & \\
$(33;14,7;7),\quad n=14$ & \\
\quad\{0,1,2,3,5,6,8,9,13,15,18,19,22,27\},\quad\{0,1,8,10,12,18,23\};&
\cite{CCS} \\
\quad\{0,2,3,4,5,8,9,11,13,16,17,19,23,26\},\quad\{0,1,2,8,13,18,22\};&$\ast$ \\
 & \\
$(33;13,12;9),\quad n=16$ & DO \\
$(33;15,11;10),\quad n=16$ & DO \\
 & \\
$(33;16,16;15),\quad n=17$ & \\
\quad\{0,1,2,3,4,5,7,8,12,13,16,18,19,21,27,28\},& \\
\quad\{0,1,2,3,5,6,8,10,12,14,15,18,22,23,25,29\};&$\ast$ \\
\quad\{0,1,2,3,4,6,8,9,10,13,14,17,20,21,25,27\},& \\
\quad\{0,1,2,3,6,8,10,11,13,15,18,19,20,21,24,30\};& Sz \\
\quad\{0,1,2,3,5,6,8,10,12,14,15,18,19,24,25,27\},& \\
\quad\{0,1,2,4,5,6,8,11,12,13,16,18,19,21,22,26\};&$\ast$ \\
 & \\
$(34;10,7;4),\quad n=13$ & \\
\quad\{0,1,3,4,7,8,13,16,18,24\},\quad\{0,2,4,9,15,16,24\};&$\ast$ \\
\quad\{0,1,3,6,8,12,15,16,17,22\},\quad\{0,1,3,9,13,17,24\};&$\ast$ \\
\quad\{0,1,3,6,10,13,15,19,26,27\},\quad\{0,2,3,4,8,14,19\};& \cite{DZ1} \\
\end{tabular}

\begin{tabular}{ll}
\multicolumn{2}{c} {\bf Table 5 (continued)} \\
 & \\
\hline
 & \\
$(34;13,7;6),\quad n=14$\quad None & \cite{MDV} \\
 & \\
$(34;12,12;8),\quad n=16$ & \\
\quad\{0,1,2,3,6,8,12,15,16,17,25,29\},& \\
\quad\{0,1,3,6,7,9,11,14,18,19,21,25\};& \cite{DZ1} \\
\quad\{0,1,2,3,7,8,13,16,18,21,25,28\},& \\
\quad\{0,1,2,4,5,8,9,11,13,19,23,25\};&$\ast$ \\
\quad\{0,1,2,4,6,11,12,15,17,20,24,28\},& \\
\quad\{0,1,4,5,6,8,11,13,14,20,21,23\};& \cite{HCD}  \\
 & \\
$(34;16,10;10),\quad n=16$ & \\
\quad\{0,1,2,3,4,5,7,8,10,14,15,19,20,25,27,28\},& \\
\quad\{0,2,4,5,8,12,16,18,21,27\};& \cite{DZ1} \\
\quad\{0,1,2,3,4,7,8,10,12,13,16,19,20,21,24,26\},& \\
\quad\{0,1,3,6,8,10,14,15,21,25\};&$\ast$ \\
\quad\{0,1,2,3,5,6,10,11,12,14,16,17,18,20,24,27\},& \\
\quad\{0,1,3,5,8,12,15,20,21,26\};&$\ast$ \\
 & \\
$(34;16,13;12),\quad n=17$ & \\
\quad\{0,1,2,3,5,6,8,12,13,14,15,18,20,22,24,31\},& \\
\quad\{0,1,4,5,7,8,9,14,15,18,23,26,28\};& \cite{DZ1} \\
\quad\{0,1,2,3,5,7,9,11,13,14,16,17,20,23,24,29\},& \\
\quad\{0,1,3,4,5,8,9,10,15,18,23,25,26\};&$\ast$ \\
\quad\{0,1,3,4,5,7,8,10,13,14,16,18,22,23,24,30\},& \\
\quad\{0,2,3,4,6,7,11,12,13,19,22,24,27\};& \cite{DZ2} \\
 & \\
$(35;8,4;2),\quad n=10$\quad None & \cite{MDV} \\
 & \\
$(35;10,4;3),\quad n=11$ & \\
\quad\{0,1,2,6,11,13,15,18,21,29\},\quad\{0,3,4,13\};& \cite{DZ1} \\
 & \\
$(35;9,6;3),\quad n=12$ & \\
\quad\{0,1,2,4,10,16,17,22,27\},\quad\{0,2,5,9,13,16\};& \cite{DZ1} \\
\quad\{0,1,2,6,8,11,14,21,25\},\quad\{0,2,5,9,17,18\};&$\ast$ \\
 & \\
$(35;12,9;6),\quad n=15$ & \\
\quad\{0,1,2,3,6,10,12,15,16,17,20,28\},\quad\{0,1,6,8,12,14,17,21,24\};&$\ast$ \\
\quad\{0,1,2,5,7,10,11,15,16,18,22,28\},\quad\{0,2,3,4,11,13,16,19,23\};&
\cite{CCS} \\
\end{tabular}

\begin{tabular}{ll}
\multicolumn{2}{c} {\bf Table 5 (continued)} \\
 & \\
\hline
 & \\
$(35;14,8;7),\quad n=15$ & \\
\quad\{0,1,2,3,4,5,10,11,13,16,18,22,26,30\},& \\
\quad\{0,1,5,7,12,16,19,22\};&$\ast$ \\
\quad\{0,1,2,4,5,7,10,12,13,17,21,23,27,28\},& \\
\quad\{0,1,4,5,7,14,16,22\};& \cite{DZ1} \\
 & \\
$(35;14,10;8),\quad n=16$ & \\
\quad\{0,1,2,3,6,8,10,12,13,16,19,20,21,25\},& \\ 
\quad\{0,1,3,4,8,13,15,21,24,29\};&$\ast$ \\
\quad\{0,1,2,4,7,9,10,12,13,17,19,23,24,28\},& \\ 
\quad\{0,1,2,3,5,8,14,18,22,28\};& \cite{DZ1} \\
 & \\
$(35;17,17;16),\quad n=18$ & \\
\quad\{0,1,2,3,4,5,8,9,10,13,14,17,19,21,24,25,27\},& \\
\quad\{0,2,3,6,7,8,9,12,13,15,17,19,20,22,26,27,29\};&$\ast$ \\
\quad\{0,1,2,3,4,6,8,9,10,11,12,17,20,23,24,27,30\},& \\
\quad\{0,2,3,4,6,7,8,12,13,15,17,20,22,25,26,27,31\};& Sz \\
 & \\
$(36;11,6;4),\quad n=13$ & \\
\quad\{0,1,2,6,8,11,13,20,21,24,28\},\quad\{0,2,3,6,12,17\};&$\ast$ \\
\quad\{0,1,3,4,6,8,12,16,17,22,29\},\quad\{0,1,7,10,19,21\};& \cite{DZ1} \\
 & \\
$(36;16,11;10),\quad n=17$ & \\
\quad\{0,1,2,3,5,6,7,10,11,13,14,21,23,24,27,29\},& \\
\quad\{0,2,3,7,9,11,14,18,19,24,30\};& \cite{MDV} \\
\quad\{0,1,2,3,6,7,8,11,12,14,17,20,21,24,26,28\},& \\
\quad\{0,1,2,6,8,9,11,13,21,24,28\};&$\ast$ \\
 & \\
$(36;15,15;12),\quad n=18$\quad None & \cite{AX} \\
 & \\
$(37;6,3;1),\quad n=8$ & \\
\quad\{0,3,5,9,17,24\},\quad\{0,1,11\};&$\ast$ \\
\quad\{0,4,6,13,18,21\},\quad\{0,1,11\};& \cite{CCS} \\
 & \\
$(37;7,6;2),\quad n=11$ & \\
\quad\{0,1,2,7,11,19,23\},\quad\{0,2,5,8,15,28\};&$\ast$ \\
\quad\{0,1,4,9,11,17,23\},\quad\{0,2,3,7,18,27\};& \cite{CCS} \\
\end{tabular}

\begin{tabular}{ll}
\multicolumn{2}{c} {\bf Table 5 (continued)} \\
 & \\
\hline
 & \\
$(37;12,4;4),\quad n=12$ & \\
\quad\{0,1,3,5,6,7,11,17,20,22,29,30\},\quad\{0,4,7,16\};& \cite{DZ1} \\
 & \\
$(37;9,9;4),\quad n=14$ & \\
\quad\{0,1,2,3,5,10,16,20,27\},\quad\{0,1,5,9,11,14,17,23,30\};& \cite{DZ1} \\
\quad\{0,1,2,6,10,15,18,20,31\},\quad\{0,2,3,4,9,12,15,19,26\};&$\ast$ \\
 & \\
$(37;10,10;5),\quad n=15$ & \\
\quad\{0,1,2,6,8,10,15,20,23,26\},\quad\{0,2,3,4,9,12,15,16,20,30\};&
\cite{CCS} \\
\quad\{0,1,4,5,11,13,14,16,20,30\},\quad\{0,1,4,6,10,15,17,18,23,25\};&
\cite{KT} \\
 & \\
$(37;15,3;6),\quad n=12$ & \\
\quad\{0,1,2,4,5,7,8,11,13,15,18,23,24,25,33\},\quad\{0,3,7\};&$\ast$ \\
\quad\{0,1,3,5,6,7,8,11,13,17,20,22,26,29,30\},\quad\{0,1,11\};& \cite{CCS} \\
 & \\
$(37;16,4;7),\quad n=13$ & \\
\quad\{0,1,2,3,6,7,9,11,12,14,16,20,21,24,27,28\},\quad\{0,2,8,22\};&$\ast$ \\
\quad\{0,2,3,4,5,7,10,11,14,16,19,23,24,25,29,31\},\quad\{0,1,4,11\};&
\cite{CCS} \\
 & \\
$(37;15,7;7),\quad n=15$ & \\
\quad\{0,1,2,5,6,8,9,11,16,18,20,21,25,29,31\},& \\
\quad\{0,2,3,8,14,15,18\};&$\ast$ \\
\quad\{0,2,3,5,6,7,11,12,16,18,19,20,22,27,30\},& \\
\quad\{0,1,4,9,11,17,23\};& \cite{CCS} \\
 & \\
$(37;13,12;8),\quad n=17$ & \\
\quad\{0,1,2,3,4,9,10,12,15,16,20,26,30\},& \\
\quad\{0,2,4,6,7,9,12,17,21,24,25,30\};& \cite{CCS} \\
\quad\{0,2,3,4,7,8,9,15,16,18,21,25,29\},& \\
\quad\{0,1,2,5,6,8,11,15,18,20,26,28\};&$\ast$ \\
 & \\
$(37;18,10;11),\quad n=17$ & \\
\quad\{0,1,2,4,5,6,7,9,10,11,16,18,19,22,26,29,30,32\},& \\
\quad\{0,1,2,6,8,10,15,20,23,26\};& \cite{CCS} \\
\quad\{0,1,2,4,5,6,7,9,11,13,16,17,19,22,23,28,29,30\},& \\
\quad\{0,2,3,5,6,10,14,19,22,29\};&$\ast$ \\
\end{tabular}

\begin{tabular}{ll}
\multicolumn{2}{c} {\bf Table 5 (continued)} \\
 & \\
\hline
 & \\
$(37;16,13;11)),\quad n=18$ & \\
\quad\{0,1,2,3,4,7,8,11,13,15,16,18,23,24,27,33\},& \\
\quad\{0,1,2,4,8,10,13,14,18,20,21,23,32\};& \cite{JC} \\
\quad\{0,1,2,3,4,7,9,13,14,17,18,19,22,25,30,31\},& \\
\quad\{0,1,3,4,6,7,9,11,16,18,20,26,30\};&$\ast$ \\
\quad\{0,1,2,3,4,8,11,12,14,17,18,20,23,25,27,32\},& \\
\quad\{0,1,3,4,5,7,11,12,16,17,20,22,30\};&$\ast$ \\
\quad\{0,1,2,4,5,7,10,12,16,18,20,23,24,25,29,32\},& \\
\quad\{0,1,2,4,5,8,11,12,13,15,21,22,27\};& \cite{MG} \\
 & \\
$(37;18,18;17),\quad n=19$ & \\
\quad\{0,1,2,3,4,5,7,8,10,14,15,16,19,21,23,28,29,32\},& \\
\quad\{0,1,2,4,5,7,8,10,11,12,16,17,20,22,24,27,28,30\};&$\ast$ \\
\quad\{0,1,2,3,4,5,7,8,11,12,14,17,18,21,28,29,30,33\},& \\
\quad\{0,1,2,3,5,8,10,14,16,18,19,20,22,24,25,27,32,34\};&$\ast$ \\
\quad\{0,1,2,3,4,6,7,9,12,14,17,18,19,20,24,27,28,31\},& \\
\quad\{0,1,2,4,5,6,8,9,10,13,15,16,20,22,24,25,31,33\};&$\ast$ \\
\quad\{0,1,2,3,5,8,9,11,13,15,16,19,21,22,23,24,28,33\},& \\
\quad\{0,1,2,4,5,6,7,9,10,11,16,18,19,22,26,29,30,32\};& \cite{JS} \\
 & \\
$(38;9,2;2),\quad n=9$\quad None & \cite{MDV} \\
 & \\
$(38;15,4;6),\quad n=13$ & \\
\quad\{0,1,2,3,5,6,9,14,15,18,20,22,25,30,32\},\quad\{0,1,7,11\};&$\ast$ \\
\quad\{0,1,2,4,5,7,10,13,15,17,21,24,25,26,33\},\quad\{0,1,7,11\};&
\cite{DZ1} \\
 & \\
$(38;12,10;6),\quad n=16$ & \\
\quad\{0,1,2,4,5,8,13,14,19,21,24,31\},& \\
\quad\{0,2,4,6,7,13,16,17,22,30\};& \cite{DZ1} \\
\quad\{0,1,4,5,7,11,13,16,18,20,21,26\},& \\
\quad\{0,1,3,4,10,14,15,22,24,30\};&$\ast$ \\
 & \\
$(38;16,8;8),\quad n=16$ & \\
\quad\{0,1,2,4,5,6,7,10,14,15,17,21,24,26,30,32\},& \\
\quad\{0,1,2,7,11,14,19,22\};&$\ast$ \\
\quad\{0,1,2,4,5,7,10,11,13,15,17,21,24,25,26,33\},& \\
\quad\{0,1,2,7,11,14,19,22\};& \cite{DZ1} \\
 & \\
$(39;8,5;2),\quad n=11$ & \\
\quad\{0,2,4,5,10,16,19,26\},\quad\{0,1,8,12,21\};&$\ast$ \\
\quad\{0,3,4,5,12,18,22,28\},\quad\{0,7,9,12,20\};& \cite{CCS} \\
\end{tabular}

\begin{tabular}{ll}
\multicolumn{2}{c} {\bf Table 5 (continued)} \\
 & \\
\hline
 & \\
$(39;9,7;3),\quad n=13$ & \\
\quad\{0,1,2,5,13,16,22,26,32\},\quad\{0,5,6,8,10,17,24\};& \cite{CCS} \\
\quad\{0,1,4,5,10,12,19,27,30\},\quad\{0,1,3,5,11,17,24\};&$\ast$ \\
 & \\
$(39;12,5;4),\quad n=13$ & \\
\quad\{0,1,3,5,9,10,16,17,18,20,23,28\},\quad\{0,4,10,13,25\};& \cite{MDV} \\
\quad\{0,2,3,4,9,10,13,15,17,22,25,31\},\quad\{0,1,4,15,20\};&$\ast$ \\
 & \\
$(39;11,7;4),\quad n=14$ & \\
\quad\{0,1,2,4,7,8,14,18,23,28,31\},\quad\{0,2,7,8,11,20,22\};&$\ast$ \\
\quad\{0,1,4,6,8,13,15,18,21,31,32\},\quad\{0,1,2,5,11,17,21\};& \cite{MDV} \\
 & \\
$(39;13,9;6),\quad n=16$ & \\
\quad\{0,1,3,6,7,9,14,18,19,21,23,29,33\},& \\
\quad\{0,1,2,4,10,11,15,18,23\};&$\ast$ \\
\quad\{0,1,4,5,7,10,11,16,18,19,24,26,28\},& \\
\quad\{0,2,3,4,8,13,16,23,30\};& \cite{CCS} \\
 & \\
$(39;15,8;7),\quad n=16$ & \\
\quad\{0,1,2,3,4,7,10,12,15,16,18,22,23,28,32\},& \\
\quad\{0,2,5,6,10,17,19,26\};&$\ast$ {} \\
\quad\{0,2,3,6,7,8,10,11,16,17,18,22,25,27,29\},& \\
\quad\{0,1,3,9,13,16,22,27\};& \cite{CCS} \\
 & \\
$(39;13,11;7),\quad n=17$ & \\
\quad\{0,1,3,5,7,9,12,16,17,20,22,23,30\},& \\
\quad\{0,1,4,5,7,12,13,15,24,25,30\};&$\ast$ \\
\quad\{0,1,4,5,7,10,11,16,18,19,24,26,28\},& \\
\quad\{0,1,2,3,5,9,13,16,22,27,32\};& \cite{CCS} \\
  & \\
$(39;15,12;9),\quad n=18$ & \\
\quad\{0,1,2,3,5,6,10,11,13,16,17,20,24,26,32\},& \\
\quad\{0,1,4,6,8,12,13,20,22,23,25,34\};&$\ast$ \\
\quad\{0,2,3,4,6,8,9,13,16,18,19,22,23,30,31\},& \\
\quad\{0,1,2,4,7,8,10,15,19,21,26,31\};& \cite{CCS} \\
 & \\
$(39;19,19;18),\quad n=20$ & \\
\quad\{0,1,2,3,4,5,6,9,13,15,17,20,22,23,25,28,30,32,36\},& \\
\quad\{0,1,2,3,5,6,9,10,11,12,15,16,17,19,23,24,29,30,33\};& Sz \\
\end{tabular}

\begin{tabular}{ll}
\multicolumn{2}{c} {\bf Table 5 (continued)} \\
 & \\
\hline
 & \\
$(39;19,19;18),\quad n=20$ & \\
\quad\{0,1,2,3,5,7,8,9,13,14,16,18,19,20,22,23,30,32,33\},& \\
\quad\{0,1,2,4,5,6,7,9,12,13,16,17,19,22,24,25,27,31,35\};&$\ast$ \\
\quad\{0,1,2,3,6,7,8,9,11,15,18,19,21,22,23,25,30,31,36\},& \\
\quad\{0,1,2,3,4,5,8,9,12,14,15,17,19,22,28,30,32,34,35\};& \cite{CCS} \\
 & \\
$(40;9,3;2),\quad n=10$ & \\
\quad\{0,1,2,5,8,13,17,19,26\},\quad\{0,10,20\};& \cite{DZ1} \\
 & \\
$(40;18,3;8),\quad n=13$\quad None & \cite{MDV} \\
 & \\
$(40;16,9;8),\quad n=17$ & \\
\quad\{0,1,2,3,6,8,10,14,17,18,19,21,24,28,29,34\},& \\
\quad\{0,1,3,4,6,10,15,23,31\};&$\ast$ \\
\quad\{0,1,3,4,5,9,10,11,15,17,18,21,25,27,30,32\},& \\
\quad\{0,2,3,4,9,13,16,21,24\};&$\ast$ \\
\quad\{0,2,3,4,5,6,8,9,12,16,17,21,23,28,30,34\},& \\
\quad\{0,3,5,10,11,19,20,27,30\};& \cite{DZ1} \\
 \\
$(40;13,13;8),\quad n=18$ & \\
\quad\{0,1,2,3,4,7,10,14,18,20,23,28,35\},& \\
\quad\{0,1,2,3,6,10,12,16,17,21,29,32,34\};&$\ast$ \\
\quad\{0,1,2,4,5,12,14,16,18,19,24,25,33\},& \\
\quad\{0,3,4,6,8,10,11,15,16,21,24,30,33\};&$\ast$ \\
\quad\{0,1,2,5,6,8,10,16,17,18,21,27,28\},& \\
\quad\{0,2,3,4,7,9,12,15,18,22,26,28,35\};&$\ast$ \\
\quad\{0,2,3,5,6,7,12,14,16,23,24,27,30\},& \\
\quad\{0,2,4,5,8,10,13,14,18,19,25,32,33\};& \cite{DZ1} \\
 \\
$(40;18,16;14),\quad n=20$ & \\
\quad\{0,1,2,3,4,6,7,8,12,14,15,17,19,22,23,26,31,32\},& \\
\quad\{0,1,2,3,6,8,10,13,14,18,21,23,24,27,34,36\};&$\ast$ \\
\quad\{0,1,3,4,6,8,9,10,13,16,18,21,22,23,29,32,33,34\},& \\
\quad\{0,1,2,5,6,7,9,11,15,17,19,20,22,23,26,33\};&$\ast$ \\
\quad\{0,2,3,4,5,6,8,11,13,14,16,17,20,21,22,26,30,33\},& \\
\quad\{0,1,2,3,7,8,11,14,18,22,23,24,26,28,33,35\};& GP \\
\quad\{0,2,3,4,5,7,11,13,16,17,19,21,22,24,28,29,35,36\},& \\
\quad\{0,1,2,4,5,8,10,11,13,14,15,20,23,24,28,30\};& \cite{DZ1} \\
\end{tabular}

\section{The range $40<v\le50$}

Our most important results are given in this section.
The SDSs having the parameter sets mentioned in the
Introduction are all included in Table 6.
We list all 89 feasible parameter sets in the above
range, and in each case list (or give reference to)
all known SDSs.

We have inserted in this table the three SDSs constructed
by Morales \cite{LM} with parameters
\[ (45;11,11;5),\quad (46;10,10;4),\quad (49;9,9;3) \]
and one of the SDSs constructed by Abel \cite{RA},
with parameters (47;7,7;3).
Morales found several families in some cases but only one 
was reported in his paper. Apparently he did not check
the families for equivalence.
In each of the four cases mentioned above we have
constructed a non-equivalent SDS.

Moreover, we found an SDS with parameters $(45;22,22;21)$.
By a well known old result of Bose \cite{RB}, it gives
a BIBD with parameters mentioned in the abstract.
However, the most interesting SDS that we have constructed is probably
the very last one in the table as it gives the first example of
a pair of binary sequences of length 50 having zero periodic
autocorrelation function.
Our genetic program was running for about 7-8 days in order
to find this particular SDS.
While Golay pairs of length 50 do not exist \cite{TA,BF}, we see
that the periodic analog of them exists. The first example of this
phenomenon has been observed earlier for length 34 (see \cite{DZ2}).

There are five D-optimal cases in this range, two with $v=43$
and one for each of $v=41,45,49$.
For $v=45$ the equivalence classes have been
enumerated in \cite{KK2}. As there are 1358 equivalence classes
of SDSs with these parameters, we refer the interested
reader to that paper.
For $v=41$ the first D-optimal SDS has been constructed in
\cite{CCS}, and we have now constructed a new one.
For $v=43$ we list the 5 D-optimal SDSs with $\la=13$
and the 10 with $\la=15$ constructed in \cite{CCS}.
We have contributed one new SDS in the former and two in the
latter case. Finally, the first D-optimal SDS with $v=49$
has been constructed in \cite{JC}, and we have found a new one.

\newpage
\begin{tabular}{ll}
\multicolumn{2}{c}
{\bf Table 6: $(v;r,s;\la)$ difference families with $40<v\le50$} \\
 & \\
\hline
 & \\
$(41;5,5;1),\quad n=9$ & \\
\quad\{0,1,4,11,29\},\quad\{0,2,8,17,22\};& \cite{KT} \\
& \\
$(41;10,6;3),\quad n=13$ & \\
\quad\{0,1,2,6,13,15,17,22,25,35\},\quad\{0,1,4,10,15,18\};&$\ast$ \\
\quad\{0,1,3,6,8,14,15,24,25,29\},\quad\{0,2,6,10,13,22\};& \cite{MDV} \\
& \\
$(41;11,10;5),\quad n=16$ & \\
\quad\{0,1,2,3,5,9,13,18,19,23,34\},\quad\{0,2,5,8,11,15,17,22,29,30\};&
\cite{CCS} \\
\quad\{0,1,2,4,5,9,12,18,20,25,31\},\quad\{0,2,3,6,8,12,19,20,29,34\};&$\ast$ \\
& \\
$(41;15,6;6),\quad n=15$ & ? \\
& \\
$(41;15,11;8),\quad n=18$ & \\
\quad\{0,1,2,3,5,7,8,9,13,18,22,23,26,32,34\},& \\
\quad\{0,2,5,8,11,15,17,18,22,29,30\};& \cite{CCS} \\
\quad\{0,1,3,4,5,7,10,11,14,15,23,26,28,32,34\},& \\
\quad\{0,1,2,4,7,9,16,17,22,28,33\};&$\ast$ \\
& \\
$(41;20,5;10),\quad n=15$ & \\
\quad\{0,1,2,3,5,6,7,9,11,13,14,16,17,21,25,26,28,31,34,35\},& \\
\quad\{0,1,6,17,19\};& \cite{MDV} \\
& \\
$(41;16,16;12),\quad n=20$ & \\
\quad\{0,1,2,3,5,7,8,9,13,18,19,22,23,26,32,34\},& \\
\quad\{0,1,3,4,6,8,11,13,15,16,17,23,24,27,30,36\};& \cite{MG} \\
\quad\{0,1,2,4,5,6,8,10,11,14,20,21,27,29,32,34\},& \\
\quad\{0,1,3,5,6,7,12,13,15,17,20,23,24,28,31,32\};&$\ast$ \\
\quad\{0,1,3,4,5,6,10,14,15,17,21,23,24,29,31,36\},& \\
\quad\{0,2,3,4,6,7,11,12,13,15,19,20,23,26,29,31\};& \cite{JC} \\
& \\
$(41;20,20;19),\quad n=21$ & \\
\quad\{0,1,2,3,4,5,7,9,13,14,16,20,21,22,26,27,30,31,33,36\},& \\
\quad\{0,1,2,3,4,6,7,9,10,11,13,14,17,19,22,25,27,29,33,34\};&$\ast$ \\
\quad\{0,1,2,3,4,6,7,9,10,11,14,18,19,20,22,23,27,30,36,37\},& \\
\quad\{0,1,2,3,5,7,9,11,13,14,15,16,21,24,26,27,30,31,33,36\};& Sz \\
\quad\{0,1,2,3,4,6,8,11,13,15,16,17,18,19,23,24,27,33,36,37\},& \\
\quad\{0,1,2,5,6,8,9,11,12,14,15,18,19,20,26,28,30,31,33,35\};& \cite{JS} \\
\end{tabular}

\newpage
\begin{tabular}{ll}
\multicolumn{2}{c} {\bf Table 6 (continued)} \\
 & \\
\hline
 & \\
$(42;16,3;6),\quad n=13$ & \\
\quad\{0,1,3,4,5,6,9,10,12,16,20,22,25,27,34,35\},\quad\{0,14,28\};&
\cite{MDV} \\
& \\
$(42;17,8;8),\quad n=17$ & \\
\quad\{0,1,2,4,7,9,10,11,13,16,17,20,21,26,28,30,31\},& \\
\quad\{0,1,5,8,13,19,25,27\};& \cite{MDV} \\
\quad\{0,1,3,4,5,7,10,12,13,14,16,22,23,27,28,30,35\},& \\
\quad\{0,2,6,10,11,16,24,27\};&$\ast$ \\
\quad\{0,2,3,4,5,6,8,12,13,14,17,19,22,26,27,29,35\},& \\
\quad\{0,1,4,8,15,20,26,32\};&$\ast$ \\
 & \\
$(42;13,10;6),\quad n=17$ & \\
\quad\{0,1,2,3,4,7,8,13,18,21,26,33,35\},& \\
\quad\{0,1,4,7,9,16,20,22,28,32\};&$\ast$ \\
\quad\{0,1,6,7,8,11,15,17,19,20,26,29,32\},& \\
\quad\{0,2,3,5,7,13,17,21,22,29\};& \cite{MDV} \\
 & \\
$(42;20,6;10),\quad n=16$ & \\
\quad\{0,1,2,3,5,6,8,10,12,13,16,18,19,20,22,25,26,31,33,34\},& \\
\quad\{0,4,5,9,20,27\};&$\ast$ \\
\quad\{0,1,2,3,5,7,9,11,12,16,17,19,20,23,24,25,27,30,32,33\},& \\
\quad\{0,3,4,9,15,32\};&$\ast$ \\
\quad\{0,1,2,4,5,6,9,12,14,15,17,19,20,21,22,26,28,29,32,38\},& \\
\quad\{0,1,8,12,19,21\};&$\ast$ \\
\quad\{0,1,2,5,6,7,8,9,11,13,17,18,20,21,23,25,28,31,34,35\},& \\
\quad\{0,1,6,10,19,21\};& \cite{MDV} \\
 & \\
$(42;21,9;12),\quad n=18$ & \\
\quad\{0,1,2,4,5,6,9,12,13,15,16,17,18,22,23,24,25,27,30,32,38\},& \\
\quad\{0,2,3,6,11,13,20,26,30\};&$\ast$ \\
\quad\{0,1,2,4,6,7,9,10,11,13,14,17,20,21,22,25,27,29,30,31,37\},& \\
\quad\{0,2,3,8,12,17,18,20,31\};& \cite{DZ1} \\
\quad\{0,1,3,4,5,8,9,10,11,13,14,16,17,21,22,24,28,31,33,35,37\},& \\
\quad\{0,1,2,4,12,17,18,24,27\};&$\ast$ \\
 & \\
$(43;6,4;1),\quad n=9$\quad None & \cite{MDV} \\
$(43;9,4;2),\quad n=11$ & ? \\
\end{tabular}

\newpage
\begin{tabular}{ll}
\multicolumn{2}{c} {\bf Table 6 (continued)} \\
 & \\
\hline
& \\
$(43;7,7;2),\quad n=12$ & \\
\quad\{0,1,3,8,12,18,24\},\quad\{0,1,5,8,19,21,34\};& \cite{RA} \\
\quad\{0,1,4,10,12,20,25\},\quad\{0,2,3,7,14,20,29\};&$\ast$ \\
& \\
$(43;13,4;4),\quad n=13$ & \\
\quad\{0,1,2,5,6,9,12,16,18,21,23,29,31\},\quad\{0,1,9,19\};& \cite{MDV} \\
& \\
$(43;16,4;6),\quad n=14$ & \\
\quad\{0,2,3,4,6,7,10,11,13,15,20,23,28,29,34,35\},& \\
\quad\{0,2,12,29\};& \cite{CCS} \\
& \\
$(43;15,7;6),\quad n=16$ & \\
\quad\{0,1,3,4,6,8,12,14,17,18,21,26,27,28,33\},& \\
\quad\{0,3,5,12,13,20,24\};&$\ast$ \\
\quad\{0,1,3,5,6,7,8,11,15,17,20,26,27,30,38\},& \\
\quad\{0,1,8,15,17,26,30\};& \cite{CCS} \\
& \\
$(43;18,6;8),\quad n=16$ & \\
\quad\{0,1,2,3,5,6,9,10,11,12,17,20,23,24,28,30,33,37\},& \\
\quad\{0,2,4,12,17,31\};& \cite{CCS} \\
\quad\{0,1,3,5,6,7,11,12,13,15,17,20,24,25,28,31,33,34\},& \\
\quad\{0,3,4,11,18,20\};&$\ast$ \\
& \\
$(43;18,9;9),\quad n=18$ & \\
\quad\{0,1,2,4,5,6,8,10,11,14,18,19,21,26,27,33,34,38\},& \\
\quad\{0,2,3,5,14,17,21,23,35\};& \cite{CCS} \\
\quad\{0,2,3,4,5,7,9,10,13,15,18,19,22,25,26,33,34,36\},& \\
\quad\{0,1,2,6,8,14,19,24,28\};&$\ast$ \\
& \\
$(43;15,15;10),\quad n=20$ & \\
\quad\{0,1,2,3,4,6,12,16,18,21,23,26,30,34,37\},& \\
\quad\{0,1,2,4,7,8,9,13,14,17,18,24,26,29,37\};& \cite{CCS} \\
\quad\{0,1,2,3,5,7,9,15,18,21,23,24,28,34,35\},& \\
\quad\{0,1,3,5,7,8,10,13,14,18,19,22,26,33,34\};&$\ast$ \\
\quad\{0,1,2,3,8,9,12,13,15,17,22,24,26,32,40\},& \\
\quad\{0,1,3,5,6,10,11,13,18,19,22,26,28,29,32\};& \cite{HCD} \\
\quad\{0,1,2,4,5,9,10,12,14,16,17,24,27,28,34\},& \\
\quad\{0,1,2,4,6,8,13,14,17,22,23,25,28,31,38\};&$\ast$ \\
\end{tabular}

\newpage
\begin{tabular}{ll}
\multicolumn{2}{c} {\bf Table 6 (continued)} \\
 & \\
\hline
& \\
$(43;21,7;11),\quad n=17$ & \\
\quad\{0,1,3,4,5,6,7,10,11,13,15,18,19,21,23,25,26,30,32,34,35\},& \\
\quad\{0,1,7,10,17,18,23\};& \cite{CCS} \\
\quad\{0,2,3,4,5,6,7,10,11,12,15,17,18,20,24,26,27,31,32,35,36\},& \\
\quad\{0,2,4,10,13,20,26\};&$\ast$ \\
& \\
$(43;18,13;11),\quad n=20$ & \\
\quad\{0,1,2,3,5,8,9,13,14,15,18,21,23,24,28,32,35,39\},& \\
\quad\{0,1,2,3,6,8,9,12,17,19,27,29,31\};&$\ast$ \\
\quad\{0,1,2,4,6,8,9,10,11,13,16,19,21,22,28,29,33,35\},& \\
\quad\{0,1,3,6,7,11,16,17,20,24,28,29,31\};&$\ast$ \\
\quad\{0,1,4,5,6,7,10,12,13,14,16,20,24,27,28,29,31,38\},& \\
\quad\{0,1,3,8,9,11,14,17,19,21,26,30,31\};& \cite{CCS} \\
& \\
$(43;18,16;13),\quad n=21$ & \\
\quad\{0,1,2,3,4,7,9,11,12,13,16,19,22,24,25,29,30,36\},& \\
\quad\{0,1,2,4,5,6,9,14,16,17,20,24,26,31,33,39\};& \cite{CK} \\
\quad\{0,1,2,3,4,7,9,11,12,13,16,19,22,24,25,29,30,36\},& \\
\quad\{0,1,3,5,8,9,14,17,18,19,21,24,25,29,31,33\};& \cite{CK} \\
\quad\{0,1,2,3,6,7,9,10,13,16,18,20,21,22,27,28,30,32\},& \\
\quad\{0,2,3,4,7,8,10,12,15,20,21,23,27,32,36,37\};& \cite{CK} \\
\quad\{0,1,2,4,5,6,9,10,11,16,17,18,21,24,26,30,34,37\},& \\
\quad\{0,1,2,3,4,7,10,12,14,16,21,22,25,28,30,33\};& \cite{CK} \\
\quad\{0,1,2,4,5,7,8,10,11,14,16,18,22,23,28,31,33,35\},& \\
\quad\{0,1,4,5,6,7,11,14,18,19,20,22,27,29,30,38\};&$\ast$ \\
\quad\{0,1,2,4,5,7,8,11,12,16,18,20,21,26,27,31,33,35\},& \\
\quad\{0,1,3,5,6,7,8,13,14,17,21,23,24,26,35,38\};& \cite{MG} \\
\quad\{0,2,4,5,6,7,10,11,14,15,19,21,22,24,30,31,35,37\},& \\
\quad\{0,1,2,3,4,7,11,13,16,17,19,22,24,25,29,36\};& \cite{CK} \\
& \\
$(43;21,15;15),\quad n=21$ & \\
\quad\{0,1,2,3,4,5,6,7,11,12,13,14,17,20,24,25,28,30,31,34,39\},& \\
\quad\{0,2,3,4,7,9,12,14,16,22,24,30,31,34,39\};&$\ast$ \\
\quad\{0,1,2,3,4,6,7,9,10,11,13,14,17,20,22,23,27,28,29,34,35\},& \\
\quad\{0,1,2,4,6,9,13,15,17,20,21,25,30,33,35\};& \cite{CK} \\
\quad\{0,1,2,3,4,6,7,9,10,11,13,14,17,20,22,23,27,28,29,34,35\},& \\
\quad\{0,1,3,5,9,11,14,15,16,19,23,28,31,33,35\};& \cite{CK} \\
\end{tabular}

\newpage
\begin{tabular}{ll}
\multicolumn{2}{c} {\bf Table 6 (continued)} \\
 & \\
\hline
& \\
$(43;21,15;15),\quad n=21$ & \\
\quad\{0,1,2,3,4,6,8,9,13,15,16,18,19,21,24,25,28,33,35,38,39\},& \\
\quad\{0,1,2,4,5,12,13,15,16,17,21,23,30,37,39\};& \cite{CK} \\
\quad\{0,1,2,3,4,6,8,9,13,15,16,18,19,21,24,25,28,33,35,38,39\},& \\
\quad\{0,1,2,5,9,10,16,17,18,20,21,23,27,29,31\};& \cite{CK} \\
\quad\{0,1,2,3,4,7,10,11,12,14,16,19,21,22,24,27,28,30,32,33,37\},& \\
\quad\{0,1,2,4,5,6,8,12,15,19,20,21,27,34,39\};& \cite{CK} \\
\quad\{0,1,2,3,4,7,10,11,12,14,16,19,21,22,24,27,28,30,32,33,37\},& \\
\quad\{0,1,5,6,7,10,12,13,14,16,20,24,28,29,31\};& \cite{CK} \\
\quad\{0,1,2,3,4,7,10,11,13,14,16,19,21,22,23,24,27,28,32,37,39\},& \\
\quad\{0,1,2,3,4,7,8,12,14,16,20,26,29,31,36\};& \cite{CK} \\
\quad\{0,1,2,3,5,6,7,9,12,15,16,18,19,20,22,26,27,28,29,34,38\},& \\
\quad\{0,1,2,5,8,10,11,13,15,19,23,24,26,31,39\};& \cite{CK} \\
\quad\{0,1,2,3,5,6,7,10,11,14,17,18,20,23,25,26,27,28,33,36,39\},& \\
\quad\{0,1,2,3,5,7,8,12,14,16,20,22,31,32,36\};& \cite{CK} \\
\quad\{0,1,2,3,5,7,10,11,12,13,14,17,19,22,25,26,27,30,33,35,39\},& \\
\quad\{0,1,3,5,6,7,9,15,16,19,20,22,26,27,34\};& \cite{CK} \\
\quad\{0,1,2,3,6,7,8,10,12,14,16,17,21,22,24,25,27,30,33,34,36\},& \\
\quad\{0,1,3,4,5,8,11,12,16,17,18,23,26,28,30\};&$\ast$ \\
\quad\{0,1,2,4,5,6,9,10,11,14,16,17,18,21,24,26,30,32,34,35,37\},& \\
\quad\{0,2,3,5,6,7,9,12,15,16,20,26,27,28,34\};& \cite{CK} \\
\quad\{0,1,2,4,5,7,8,9,10,13,14,17,18,19,22,24,26,29,30,32,37\},& \\
\quad\{0,1,2,3,4,6,11,13,17,20,23,27,29,34,35\};& \cite{CK} \\
\quad\{0,1,3,4,5,6,8,9,10,14,15,17,18,21,25,27,28,32,33,35,37\},& \\
\quad\{0,2,3,4,6,8,9,13,15,16,21,26,28,29,37\};& \cite{JC} \\
& \\
$(43;21,21;20),\quad n=22$ & \\
\quad\{0,1,2,3,4,6,8,9,11,14,15,18,20,23,24,25,28,33,34,36,40\},& \\
\quad\{0,1,2,3,5,7,8,9,11,13,14,15,18,19,22,26,28,29,30,31,38\};&$\ast$ \\
\quad\{0,1,2,3,5,6,9,10,11,12,17,20,21,23,24,26,28,30,33,34,37\},& \\
\quad\{0,1,2,4,5,6,7,8,9,13,14,15,18,20,21,23,28,31,33,35,39\};& \cite{DZ1} \\
& \\
$(44;8,6;2),\quad n=12$\quad None & \cite{MDV} \\
$(44;19,2;8),\quad n=13$ & ? \\
\end{tabular}

\newpage
\begin{tabular}{ll}
\multicolumn{2}{c} {\bf Table 6 (continued)} \\
 & \\
\hline
& \\
$(44;17,9;8)),\quad n=18$ & \\
\quad\{0,1,2,3,4,7,9,11,15,17,20,22,23,27,28,32,37\},& \\
\quad\{0,1,4,10,12,13,16,23,30\};&$\ast$ \\
\quad\{0,1,2,4,6,7,12,13,18,20,22,23,26,27,30,35,37\},& \\
\quad\{0,2,3,6,10,15,16,18,27\};& \cite{MDV} \\
\quad\{0,1,2,5,6,7,8,9,15,18,19,22,24,26,29,34,38\},& \\
\quad\{0,1,4,6,12,14,17,26,35\};&$\ast$ \\
& \\
$(44;18,15;12),\quad n=21$ \quad None & \cite{MDV} \\
$(44;21,14;14),\quad n=21$ \quad None & \cite{MDV} \\
$(45;7,2;1),\quad n=8$ \quad None & \cite{MDV} \\
& \\
$(45;10,7;3),\quad n=14$ & \\
\quad\{0,1,2,4,9,11,14,19,25,31\},\quad\{0,4,8,11,17,29,30\};&$\ast$ \\
\quad\{0,3,4,5,10,17,19,21,27,30\},\quad\{0,1,6,9,13,21,35\};& \cite{MDV} \\
& \\
$(45;13,5;4),\quad n=14$ & \\
\quad\{0,1,4,6,8,12,13,14,21,24,30,31,34\},\quad\{0,3,5,14,19\};& \cite{MDV} \\
& \\
$(45;11,11;5),\quad n=17$ & \\
\quad\{0,1,2,8,10,11,14,19,24,27,31\},& \\
\quad\{0,2,3,4,7,9,16,22,27,33,37\};&$\ast$ \\
\quad\{0,1,4,7,8,16,17,22,27,31,33\},& \\
\quad\{0,1,5,7,8,10,12,18,20,29,32\};& \cite{LM} \\
& \\
$(45;12,12;6),\quad n=18$ & \\
\quad\{0,1,2,3,8,10,14,15,19,24,27,30\},& \\
\quad\{0,1,3,4,6,10,14,19,26,28,34,38\};& \cite{MDV} \\
\quad\{0,1,2,6,7,9,13,18,21,23,26,36\},& \\
\quad\{0,1,4,7,8,10,16,18,20,27,31,32\};&$\ast$ \\
& \\
$(45;18,2;7),\quad n=13$ & ? \\
& \\
$(45;18,10;9),\quad n=19$ & \\
\quad\{0,1,2,5,7,9,11,12,13,15,16,20,21,28,30,33,34,36\},& \\
\quad\{0,1,4,6,9,15,16,23,26,33\};&$\ast$ \\
\end{tabular}

\newpage
\begin{tabular}{ll}
\multicolumn{2}{c} {\bf Table 6 (continued)} \\
 & \\
\hline
& \\
$(45;16,13;9),\quad n=20$ & \\
\quad\{0,1,2,5,10,11,12,16,17,18,20,23,25,31,34,37\},& \\
\quad\{0,1,3,4,5,10,13,17,20,22,24,28,32\};& \cite{MDV} \\
\quad\{0,1,3,4,6,10,12,14,15,17,21,22,26,27,30,32\},& \\
\quad\{0,1,2,8,9,10,14,18,21,25,28,31,33\};&$\ast$ \\
& \\
$(45;21,5;10),\quad n=16$ & \\
\quad\{0,1,2,3,4,5,6,9,10,14,15,17,21,22,23,25,28,31,35,38,40\},& \\
\quad\{0,4,6,15,33\};& \cite{MDV} \\
& \\
$(45;22,11;13),\quad n=20$ & \\
\quad\{0,1,2,3,4,7,9,11,14,15,18,19,20,21,23,26,27,28,31,34,36,42\},& \\
\quad\{0,1,2,3,7,10,12,16,25,31,35\};& \cite{MDV} \\
\quad\{0,2,3,4,5,6,8,9,11,13,16,17,18,19,23,26,29,30,34,35,38,39\},& \\
\quad\{0,1,3,7,9,11,14,21,23,28,29\};&$\ast$ \\
& \\
$(45;21,16;15),\quad n=22$ & DO \\
& \\
$(45;22,22;21),\quad n=23$ & \\
\quad\{0,1,2,3,6,7,8,10,11,13,14,18,20,22,23,24,26,29,31,32,33,41\},& \\
\quad\{0,1,3,4,5,6,7,11,12,14,15,18,19,20,21,24,28,30,32,35,37,40\};&$\ast$ \\
 & \\
$(46;10,10;4),\quad n=16$ & \\
\quad\{0,1,2,5,10,13,19,25,29,31\},\quad\{0,1,2,7,9,12,16,20,23,33\};&
\cite{LM} \\
\quad\{0,1,3,4,13,15,19,20,26,28\},\quad\{0,3,5,8,12,13,19,23,29,37\};&$\ast$ \\
 & \\
$(46;16,6;6),\quad n=16$ & \\
\quad\{0,1,2,5,7,8,9,12,15,20,21,23,25,27,36,37\},& \\
\quad\{0,3,9,17,22,26\};& \cite{MDV} \\
 & \\
$(46;21,6;10),\quad n=17$ & ? \\
$(46;16,15;10),\quad n=21$ \quad None & \cite{MDV} \\
$(46;21,15;14),\quad n=22$ \quad None & \cite{MDV} \\
$(47;10,2;2),\quad n=10$ \quad None & \cite{MDV} \\
$(47;9,5;2),\quad n=12$ & ? \\
$(47;12,3;3),\quad n=12$ & ? \\
$(47;14,2;4),\quad n=12$ & ? \\
$(47;15,5;5),\quad n=15$ & ? \\
\end{tabular}

\newpage
\begin{tabular}{ll}
\multicolumn{2}{c} {\bf Table 6 (continued)} \\
 & \\
\hline
& \\
$(47;19,9;9),\quad n=19$ & \\
\quad\{0,1,2,3,4,6,7,10,14,16,17,19,23,24,26,28,34,36,42\},& \\
\quad\{0,4,5,9,10,18,21,29,36\};&$\ast$ \\
\quad\{0,1,2,3,6,7,9,11,13,14,17,21,22,23,29,31,32,34,36\},& \\
\quad\{0,4,5,9,12,15,21,28,38\};& \cite{MDV} \\
& \\
$(47;22,10;12),\quad n=20$ & \\
\quad\{0,1,2,3,4,5,7,9,10,12,15,17,18,19,23,24,27,28,30,34,35,38\},& \\
\quad\{0,2,8,9,13,16,22,26,28,38\};&$\ast$ \\
\quad\{0,1,2,4,5,7,8,11,12,14,17,19,22,24,25,26,31,32,33,35,39,43\},& \\
\quad\{0,1,2,5,10,11,14,16,27,29\};& \cite{MDV} \\
& \\
$(47;21,12;12),\quad n=21$ & \\
\quad\{0,1,2,3,4,6,8,9,11,12,14,16,20,21,25,26,29,32,33,36,39\},& \\
\quad\{0,1,3,5,10,14,15,16,23,29,31,41\};&$\ast$ \\
\quad\{0,1,2,4,5,6,11,12,14,16,18,20,23,24,27,31,33,34,36,39,44\},& \\
\quad\{0,1,2,5,8,10,11,12,19,25,26,31\};&$\ast$ \\
& \\
$(47;19,15;12),\quad n=22$ & \\
\quad\{0,1,2,4,6,8,11,12,15,18,19,20,22,25,27,28,30,35,36\},& \\
\quad\{0,2,3,5,8,9,12,14,18,22,23,27,33,34,35\};&$\ast$ \\
& \\
$(47;22,14;14),\quad n=22$ & \\
\quad\{0,1,2,3,4,6,8,9,11,12,14,15,16,20,21,22,27,30,31,34,36,38\},& \\
\quad\{0,1,2,5,8,10,15,19,23,25,26,31,35,38\};&$\ast$ \\
& \\
$(47;23,23;22),\quad n=24$ & DS \\
$(48;14,3;4),\quad n=13$ & ? \\
& \\
$(48;12,8;4),\quad n=16$ & \\
\quad\{0,1,3,4,9,11,16,17,21,27,31,40\},\quad\{0,2,3,7,9,22,25,36\};&
\cite{MDV} \\
\quad\{0,2,4,5,10,12,13,18,22,27,33,34\},\quad\{0,3,4,7,14,16,27,33\};&$\ast$ \\
& \\
$(48;16,7;6),\quad n=17$ & \\
\quad\{0,1,2,3,7,9,10,14,19,20,22,25,28,32,34,36\},& \\
\quad\{0,3,4,8,13,24,41\};& \cite{MDV} \\
\quad\{0,1,3,4,8,10,11,12,15,17,23,25,28,29,33,39\},& \\
\quad\{0,2,3,9,18,31,36\};&$\ast$ \\
\end{tabular}

\newpage
\begin{tabular}{ll}
\multicolumn{2}{c} {\bf Table 6 (continued)} \\
 & \\
\hline
& \\
$(48;15,9;6),\quad n=18$ & \\
\quad\{0,1,2,4,5,10,12,14,17,20,23,24,28,31,37\},& \\
\quad\{0,1,2,6,11,18,20,26,33\};&$\ast$ \\
\quad\{0,1,2,7,9,12,14,16,22,25,30,31,34,35,38\},& \\
\quad\{0,1,3,4,9,11,15,21,32\};& \cite{MDV} \\
& \\
$(48;20,10;10),\quad n=20$ & \\
\quad\{0,1,2,6,8,9,10,11,12,15,16,19,21,24,27,29,32,34,36,38\},& \\
\quad\{0,1,4,5,11,17,18,25,34,37\};&$\ast$ \\
& \\
$(48;24,4;12),\quad n=16$ & \\
\quad\{0,1,2,4,5,6,8,9,11,13,15,16,18,19,21,24,27,28,31,32,33,& \\
\quad \quad 38,39,40\},\quad\{0,2,6,20\};&$\ast$ \\
& \\
$(49;7,3;1),\quad n=9$ & \\
\quad\{0,1,7,10,15,27,31\},\quad\{0,2,13\};& \cite{CCS} \\
& \\
$(49;10,3;2),\quad n=11$ & ? \\
& \\
$(49;12,4;3),\quad n=13$ & \\
\quad\{0,1,2,7,10,12,15,18,22,24,31,35\},\quad\{0,1,19,23\};& \cite{CCS} \\
& \\
$(49;9,9;3),\quad n=15$ & \\
\quad\{0,1,2,4,7,12,20,27,36\},\quad\{0,1,4,9,16,18,22,28,39\};& \cite{LM} \\
\quad\{0,1,3,6,12,13,20,28,38\},\quad\{0,2,4,7,8,18,23,27,36\};&$\ast$ \\
& \\
$(49;15,6;5),\quad n=16$ & \\
\quad\{0,1,2,3,4,6,10,13,17,21,22,27,29,35,40\},& \\
\quad\{0,5,12,15,21,29\};& \cite{CCS} \\
\quad\{0,1,3,7,9,10,13,15,18,22,23,32,34,38,39\},& \\
\quad\{0,4,5,7,18,26\};&$\ast$ \\
& \\
$(49;13,12;6),\quad n=19$ & \\
\quad\{0,1,2,4,6,10,11,18,22,23,30,33,36\},& \\
\quad\{0,1,3,6,11,12,15,19,21,26,28,36\};& \cite{CCS} \\
\quad\{0,1,3,5,6,9,17,19,21,25,26,36,43\},& \\
\quad\{0,1,2,7,10,12,15,21,22,25,29,38\};&$\ast$ \\
\end{tabular}

\newpage
\begin{tabular}{ll}
\multicolumn{2}{c} {\bf Table 6 (continued)} \\
 & \\
\hline
& \\
$(49;18,6;7),\quad n=17$ & \\
\quad\{0,1,2,3,7,9,12,15,16,17,20,23,27,28,30,32,34,40\},& \\
\quad\{0,5,6,9,19,28\};& \cite{MDV} \\
& \\
$(49;19,7;8),\quad n=18$ & \\
\quad\{0,1,2,3,5,7,8,9,12,16,18,19,22,24,30,31,32,36,41\},& \\
\quad\{0,3,8,12,15,26,36\};&$\ast$ \\
\quad\{0,1,2,4,5,8,10,13,15,17,19,23,24,29,30,31,34,37,41\},& \\
\quad\{0,3,5,6,14,15,31\};& \cite{MDV} \\
& \\
$(49;21,4;9),\quad n=16$ & ? \\
& \\
$(49;19,10;9),\quad n=20$ & \\
\quad\{0,1,2,3,5,8,10,12,13,16,20,21,22,25,27,33,34,36,40\},& \\
\quad\{0,4,6,7,10,14,22,23,28,33\};& \cite{MDV} \\
\quad\{0,1,2,4,5,7,9,11,16,17,19,20,23,26,27,31,32,37,40\},& \\
\quad\{0,1,3,5,11,13,18,24,25,33\};&$\ast$ \\
& \\
$(49;16,16;10),\quad n=22$ & \\
\quad\{0,1,2,4,6,9,10,12,13,18,20,23,30,34,35,43\},& \\
\quad\{0,1,3,4,5,6,11,12,18,22,25,27,31,35,38,40\};&$\ast$ \\
\quad\{0,1,3,4,7,8,12,14,15,16,25,27,30,32,35,44\},& \\
\quad\{0,1,4,6,7,8,14,16,17,19,23,25,29,30,35,39\};& \cite{MDV} \\
& \\
$(49;21,13;12),\quad n=22$ & \\
\quad\{0,1,2,3,5,7,8,9,13,14,16,17,18,23,26,28,29,32,33,36,45\},& \\
\quad\{0,1,3,6,11,13,15,17,23,30,31,38,41\};&$\ast$ \\
& \\
$(49;24,9;13),\quad n=20$ & \\
\quad\{0,1,2,3,5,6,8,9,11,12,17,18,19,21,23,25,28,31,32,33,36,& \\
\quad \quad 37,40,44\},\quad\{0,1,3,10,11,16,21,23,25\};&$\ast$ \\
\quad\{0,1,2,3,6,7,9,11,12,13,15,16,17,19,21,26,28,29,30,33,36,& \\
\quad \quad 37,39,44\},\quad\{0,3,4,5,12,15,21,29,34\};& \cite{CCS} \\
\end{tabular}

\newpage
\begin{tabular}{ll}
\multicolumn{2}{c} {\bf Table 6 (continued)} \\
 & \\
\hline
& \\
$(49;22,15;14),\quad n=23$ & \\
\quad\{0,1,2,3,4,5,7,9,10,14,15,17,19,22,24,25,28,31,33,35,39,41\},& \\
\quad \quad\{0,3,4,6,7,11,12,18,19,20,24,29,33,40\};&$\ast$ \\
\quad\{0,2,3,4,5,8,10,12,13,14,16,20,21,22,23,25,27,32,33,36,& \\
\quad \quad 38,39\},\quad\{0,1,5,6,8,9,12,15,19,24,27,29,33,34,41\};&
\cite{MDV} \\
& \\
$(49;22,18;16),\quad n=24$ & \\
\quad\{0,1,2,3,4,5,6,9,11,13,14,19,20,21,23,26,27,30,35,38,40,42\},& \\
\quad\{0,1,3,4,5,8,9,13,15,19,21,24,26,27,30,37,43,44\};& \cite{JC} \\
\quad\{0,1,2,3,4,5,7,10,12,13,15,18,19,22,23,24,29,30,34,36,38,43\},& \\
\quad\{0,1,3,4,6,7,9,10,14,16,17,24,26,28,32,36,37,41\};&$\ast$ \\
& \\
$(49;24,24;23),\quad n=25$ & \\
\quad\{0,1,2,3,4,5,8,9,11,13,15,17,18,19,20,23,26,27,28,30,33,37,& \\
\quad \quad 39,44\},\quad\{0,1,2,3,5,6,7,10,13,14,15,17,21,22,23,26,27,30,& \\
\quad \quad 32,33,35,36,38,44\};&$\ast$ \\
\quad\{0,1,2,4,5,6,8,9,14,15,16,17,18,20,23,25,26,27,30,31,33,35,& \\
\quad \quad 38,44\},\quad\{0,1,2,5,7,8,9,10,11,13,17,18,21,22,23,25,27,30,& \\
\quad \quad 32,33,36,37,39,43\};& \cite{MDV} \\
& \\
$(50;8,7;2),\quad n=13$ & ? \\
$(50;20,4;8),\quad n=16$ & ? \\
$(50;15,14;8),\quad n=21$ \quad None & \cite{MDV} \\
$(50;20,11;10),\quad n=21$ \quad None & \cite{MDV} \\
$(50;20,18;14),\quad n=24$ \quad None & \cite{MDV} \\
$(50;22,21;18),\quad n=25$ & ? \\
& \\
$(50;25,20;20),\quad n=25$ & \\
\quad\{0,1,2,3,5,6,7,9,10,11,13,15,17,20,22,23,26,27,28,29,31,36,& \\
\quad \quad 38,39,45\}, & \\
\quad\{0,1,2,3,5,6,8,14,15,17,18,21,25,27,32,35,36,40,44,45\}.&$\ast$ \\
\end{tabular}


\begin{thebibliography}{99}

\bibitem{RA}
R.J.R. Abel, Forty-three balanced incomplete block designs,
J. Combin. Theory Ser. A {\bf 65} (1994), 252--267.

\bibitem{TA}
T.H. Andres, Some combinatorial properties of complementary
sequences, M.Sc. Thesis, University of Manitoba, Winnipeg, 1977.

\bibitem{AX}
K.T. Arasu and Q. Xiang, On the existence of periodic complementary
binary sequences,
Designs, Codes and Cryptography {\bf 2} (1992), 257--262.

\bibitem{DA}
D. Ashlock, Finding designs with genetic algorithms, in W.D. Wallis
(Ed.), Computational and Constructive Design Theory, pp. 49--65,
Kluwer Academic Publishers, Dordrecht/Boston/London, 1996.

\bibitem{BA}
L. B\"omer and M. Antweiler, Periodic complementary binary sequences,
IEEE Trans. Inform. Theory {\bf 36} (1990), 1487--1494.

\bibitem{BF}
P.B. Borwein and R.A. Ferguson, A complete description of Golay
pairs for lengths up to 100,
Math. Comp. {\bf 73} (2003), 967--985.

\bibitem{RB}
R.C. Bose, On the construction of balanced incomplete block
designs, Ann. Eugenics {\bf 9} (1939), 353--399.

\bibitem{CCS}
S. Chadjiconstantinidis, T. Chadjipadelis and K. Sotirakoglou,
Two cyclic supplementary difference sets and optimal designs in
linear models,
J. Comb. Math. Comb. Comput. {\bf 18} (1995), 33--56.

\bibitem{CK}
T. Chadjipantelis and S. Kounias, Supplementary difference sets and
D-optimal designs for $n\equiv2$ mod 4,
Discrete Math. {\bf 57} (1985), 211--216.

\bibitem{JC}
J.H.E. Cohn, On determinants with elements $\pm1$, II,
Bull. London Math. Soc. {\bf 21} (1989), 36--42.

\bibitem{HCD}
C.J. Colbourn and J.H. Dinitz, Editors, Handbook of Combinatorial Designs,
	2nd edition, Chapman \& Hall, Boca Raton/London/New York, 2007.

\bibitem{CD}
C. Ding, Two constructions of $(v,(v-1)/2,(v-3)/2)$ difference families,
J. Combin. Designs {\bf 16} (2008), 164--171.

\bibitem{DZ1}
D.\v{Z}. \DJo{},
Survey of cyclic $(v;r,s;\la)$ difference families with $v\le50$,
Facta Universitatis (Ni\v{s}), Ser. Mathematics and Informatics
{\bf 12} (1997), 1--13.

\bibitem{DZ2}
\bysame, Note on periodic complementary sets of binary sequences,
Designs, Codes and Cryptography {\bf 13} (1998), 251--256.


\bibitem{DZ3}
\bysame, Equivalence classes and representatives of Golay sequences,
Discrete Math. {\bf 189} (1998), 79--93.

\bibitem{HE}
H. Ehlich, Determinantenabsch\"atzungen f\"ur bin\"are Matrizen,
Math. Zeitschr. {\bf 83} (1964), 123--132.

\bibitem{FSX}
K. Feng, P.J-S. Shiue, and Q. Xiang, On aperiodic and periodic
complementary binary sequences,
IEEE Trans. Inform. Theory {\bf 45} (1999), 296--303.

\bibitem{MG}
M. Gysin, New D-optimal designs via cyclotomy and generalised cyclotomy,
Australasian Journal of Combinatorics {\bf 15} (1997), 247--255.

\bibitem{GS1}
M. Gysin and J. Seberry, An experimental search and new combinatorial
designs via a generalisation of cyclotomy,
J. Comb. Math. Comb. Comput. {\bf 27} (1998), 143--160.

\bibitem{GS2}
M. Gysin and J. Seberry, On new families of supplementary difference
sets over rings with short orbits,
J. Comb. Math. Comb. Comput. {\bf 28} (1998), 161--186.

\bibitem{KKS}
C. Koukouvinos, S. Kounias and J. Seberry, Supplementary difference
sets and optimal designs,
Discrete Math. {\bf 88} (1991), 49--58.

\bibitem{KK1}
S. Kounias, C. Koukouvinos, N. Nicolaou and A. Kakos, The
non-equivalent circulant D-optimal designs for $n\equiv2$ mod 4,
$n\le54$, $n=66$,
J. Combin. Theory Ser. A {\bf 65} (1994), 26--38.

\bibitem{KK2}
\bysame, The non-equivalent circulant D-optimal designs for $n=90$,
J. Statist. Plann. Inference {\bf 53} (1996), 253--259.

\bibitem{MDV}
L. Mart\'{\i}nez, D.\v{Z}. \DJo{} and A. Vera-L\'{o}pez,
Existence question for difference families and construction of some
new families,
J. Combin. Designs {\bf 12} (2004), 256--270.

\bibitem{LM}
L.B. Morales, Constructing difference families through an
optimization approach: Six new BIBDs, 
J. Combin. Designs {\bf 8} (2000), 261--273.

\bibitem{AP}
A. Pott, Finite Geometry and Character Theory, Lecture Notes Math.
1601, Springer, New York, 1995.

\bibitem{SY}
J. Seberry and M. Yamada, Hadamard matrices, sequences and block
designs, in Contemporary Design Theory: A Collection of Surveys,
Eds. J.H. Dinitz and D.R. Stinson, J. Wiley, New York, 1992,
pp. 431--560.

\bibitem{JS}
J. Seberry Wallis, Some remarks on supplementary difference sets.
In Colloquia Mathematica Societatis Janos Bolyai
{\bf 10} (1973), 1503--1526.

\bibitem{KT}
K. Takeuchi, A table of difference sets generating balanced incomplete
block designs, 
Review Intern. Statistical Inst. {\bf 30} (1962), 361--366.

\bibitem{RW}
R.M. Wilson, Cyclotomy and difference families in elementary
abelian groups,
Journal of Number Theory {\bf 4} (1972), 17--47.

\bibitem{CY1}
C.H. Young, On designs of maximal $(+1,-1)$-matrices of order
$n\equiv2 \pmod{4}$,
Math. Comp. {\bf 23} (1968), 174--180.

\bibitem{CY2}
\bysame, Maximal binary matrices and sum of two squares,
Math. Comp. {\bf 30} (1976), 361--366.

\end{thebibliography}
\end{document}